\documentclass[twoside,a4paper,12pt,centertags]{amsart}
\usepackage{amsmath,amssymb,verbatim,vmargin}
\usepackage[bookmarks=true]{hyperref}   

\theoremstyle{plain}
\newtheorem{thm}{Theorem}[section]
\newtheorem{lem}[thm]{Lemma}
\newtheorem{cor}[thm]{Corollary}
\newtheorem{prop}[thm]{Proposition}

\theoremstyle{definition}
\newtheorem{defn}[thm]{Definition}
\newtheorem{rem}[thm]{Remark}

\newtheorem{ex}[thm]{Consequence}

\title[Quadratic estimates and functional calculi]{Quadratic estimates
and functional calculi of perturbed Dirac operators}
\author{Andreas Axelsson} \author{Stephen Keith} \author{Alan \Mcc Intosh}
\address{Centre for Mathematics and its Application \\
Australian National University \\ Canberra, ACT 0200, Australia}
\email{kax74@yahoo.se}
\email{stephen.keith@anu.edu.au} \email{Alan.McIntosh@anu.edu.au}

\newcommand{\qedend}{ }

\newcommand{\semic}{{:}}       
\newcommand{\Mcc}{{M\raise.55ex\hbox{\lowercase{c}}}}
\newcommand{\dirac}{{\mathbf D}}

\newcommand{\rnum}{{\mathbf R}}
\newcommand{\cnum}{{\mathbf C}}
\newcommand{\nnum}{{\mathbf N}}
\newcommand{\znum}{{\mathbf Z}}
\newcommand{\mH}{{\mathcal H}}

\newcommand{\mK}{{\mathcal K}}

\newcommand{\mL}{{\mathcal L}}

\newcommand{\mV}{{\mathcal V}}

\DeclareMathOperator{\re}{Re}
\newcommand{\im}{\text{{\rm Im}}\,}
\newcommand{\sett}[2]{ \{ #1 \, \semic \, #2 \} }

\newcommand{\brac}[1]{\langle #1 \rangle}
\newcommand{\supp}{\text{{\rm supp}}\,}
\newcommand{\dist}{\text{{\rm dist}}\,}

\newcommand{\nul}{\textsf{N}}
\newcommand{\ran}{\textsf{R}}
\newcommand{\dom}{\textsf{D}}

\newcommand{\clos}[1]{\overline{#1}}
\newcommand{\conj}[1]{\overline{#1}}
\newcommand{\dyadic}{\triangle}
\newcommand{\sgn}{\text{{\rm sgn}}}
\newcommand{\barint}{\mbox{$ave \int$}}
\newcommand{\PP}{{\mathbf P}_B}
\newcommand{\PQ}{{\mathbf P}}

\newcommand{\divv}{{\text{div}}}
\newcommand{\esssup}{\text{{\rm ess sup}}}

\def\barint_#1{\mathchoice
            {\mathop{\vrule width 6pt
height 3 pt depth -2.5pt
                    \kern -8.8pt
\intop}\nolimits_{#1}}%
            {\mathop{\vrule width 5pt height
3 pt depth -2.6pt
                    \kern -6.5pt
\intop}\nolimits_{#1}}%
            {\mathop{\vrule width 5pt height
3 pt depth -2.6pt
                    \kern -6pt
\intop}\nolimits_{#1}}%
            {\mathop{\vrule width 5pt height
3 pt depth -2.6pt
          \kern -6pt \intop}\nolimits_{#1}}}

\newcommand{\Pa}{{\tilde{\mathbf P}}_B^1}
\newcommand{\Pb}{{\tilde{\mathbf P}}_B^2}
\newcommand{\C}{\cnum}
\newcommand{\R}{\rnum}
\newcommand{\N}{\nnum}

\DeclareMathOperator{\I}{I}
\newcommand{\PPT}{T}            
\newcommand{\mC}{{\mathcal C}}
\newcommand{\CB}{{\mathcal B}}

\date{May 2004}

\begin{document}

\begin{abstract}
We prove quadratic estimates for complex perturbations of
Dirac-type operators, and thereby show that such operators have a
bounded functional calculus. As an application we show that
spectral projections of the Hodge--Dirac operator  on  compact
manifolds depend analytically on $L_\infty$ changes in the metric.
We also recover a unified proof of many results in the Calder\'on
program, including the Kato square root problem and the
boundedness of the Cauchy operator on Lipschitz curves and
surfaces.
\end{abstract}

\maketitle

\tableofcontents
\section{Introduction}

We prove quadratic estimates
\begin{equation}
\int_0^\infty\|\Pi_B(\I+t^2{\Pi_B}^2)^{-1}u\|^2\, t\, dt\ \approx\
\|u\|^2 \label{quad}
\end{equation}
for all $u \in L_2(\rnum^n, \Lambda)$, where $ \Pi_B=d+B^{-1}d^*B$
is the perturbation of a Dirac-type operator $\Pi=d+d^*$ by an
operator $B$ of multiplication  by an $L_\infty$ complex
matrix-valued function with uniformly positive real part. Here
$\Lambda$ is the complex exterior algebra on $\rnum^n$ and $d$
denotes the exterior derivative.

This estimate implies that $\Pi_B$ has a bounded functional
calculus. This means that
\begin{equation}\| f(\Pi_B)u\| \lesssim\|f\|_{_\infty}\|u\|  \label{fun}
\end{equation}
for all $u \in L_2(\rnum^n, \Lambda)$ and
    all bounded holomorphic functions
$f :  S_\mu^o \longrightarrow \cnum$, where $S_\mu^o$ is an open
double sector
$$S_\mu^o:=\{z\in\cnum : |\arg(\pm z)|<\mu\} \quad \text{with} \quad
\mu>\omega:=\sup|\arg(Bu,u)|.$$ This result in turn implies
perturbation estimates of the form
\begin{equation}
\|f(\Pi_{B+A})u-f(\Pi_B)u\|
\lesssim\|f\|_{_\infty}\|A\|_{_\infty}\|u\|
\label{morefun}\end{equation} for all $u \in L_2(\rnum^n,
\Lambda)$, provided $\| A \|_\infty$ is not too large.

The unperturbed operator $\Pi$ is selfadjoint, so when $B=I$,
(\ref{fun}) holds for all bounded Borel measurable functions $f$
by the spectral theory of selfadjoint operators. When $B$ is
positive selfadjoint, then $\Pi_B$ is selfadjoint with respect to
the inner-product $(Bu,v)$ on $L_2(\rnum^n, \Lambda)$, so
(\ref{quad}) and (\ref{fun}) still hold by spectral theory.
However (\ref{morefun}) would not, were it not for the structure
of the operators $\Pi$, $B$ and $A$. This is because  we need
(\ref{fun}) for all small non--selfadjoint perturbations of $B$ in
order to deduce (\ref{morefun}) for small selfadjoint
perturbations.

Under our assumptions on $B$, the operator $\Pi_B$ has spectrum in
the closed double sector $S_\omega=\{z\in\cnum: |\arg(\pm
z)|\leq\omega\}$ and satisfies resolvent bounds $$\|(\Pi_B-\lambda
\I)^{-1}\|\lesssim \frac{1}{\dist(\lambda,S_\omega)} $$for all
$\lambda \in \cnum \setminus S_\omega.$ This follows from operator
theory, but a proof of the quadratic estimate (\ref{quad})
requires the full strength of the harmonic analysis. Once the
estimate (\ref{quad}) is proven, then (\ref{fun}) follows if
$\omega<\mu<\frac\pi2$. It can then be seen that $f(\Pi_B)$
depends holomorphically on $B$, from which (\ref{morefun}) follows
provided $A$ is not too large.

Our result was inspired by the proof of the Kato square root
problem by Auscher, Hofmann, Lacey, \Mcc Intosh and Tchamitchian
\cite{AHLMcT}, and includes not only it as a corollary, but also
many results in the Calder\'on program such as the boundedness of
the Cauchy operator on Lipschitz curves and surfaces. The proof
uses many of the concepts developed over the years to prove these
results in the Calder\'on program, and in particular the proof of
the Kato problem, but is not a direct consequence, as the operator
$\Pi_B$ is first order, and the second order operator ${\Pi_B}^2$
is not in divergence form. Indeed, our arguments utilize only the
first order structure of the operator. This enables us to exploit
the algebra involved in the (non--orthogonal) Hodge decomposition
of the first order system
$$L_2(\rnum^n, \Lambda) =\nul(d)\oplus \nul(B^{-1}d^*B)$$ where
$\nul(d)$ is the null-space of $d$.

Combining the Hodge decomposition with (\ref{fun}) in the case
when $f(z)=z/ \sqrt{z^2}$, we obtain the equivalence $$\|du\| +
\|d^* Bu\|\approx \|\Pi_B u\|\approx \|\sqrt{{\Pi_B}^2}u\|.$$ The
square root problem of Kato  follows
    in the special case when $B$ splits as
    $B_k(x):\Lambda^k\to\Lambda^k$ for each $0\leq k\leq n$, and for
     almost every~$x\in\rnum^n$, with $B_0=I$ and
$B_1(x)=A(x):\cnum^n\to\cnum^n$. On making the identification
$$d:L_2(\rnum^n,\Lambda^0)\to L_2(\rnum^n,\Lambda^1) \quad \text{with}
\quad \nabla:L_2(\rnum^n,\cnum)\to L_2(\rnum^n,\cnum^n)\quad
\text{and}$$
$$d^*:L_2(\rnum^n,\Lambda^1)\to L_2(\rnum^n,\Lambda^0) \quad
\text{with} \quad -\divv:L_2(\rnum^n,\cnum^n)\to
L_2(\rnum^n,\cnum)$$ and restricting our attention to $\Lambda^0$,
we obtain $$\|\nabla u\|\approx\|\sqrt{-\text{div} A\nabla}u\| $$
for all $u\in L_2(\rnum^n,\cnum)$.

The choice of test--functions used in our proof of the stopping
time argument in Section \ref{harmansec} has more in common with
that presented in the paper on elliptic systems \cite{AHMcT} than
with \cite{AHLMcT}, but the result stated above does not include
the full result on systems. To remedy this, as well as to allow
further consequences, our results can in fact be stated somewhat
more generally than so far indicated, though without much effect
on the proofs. Rather than $d$, we consider any first order system
$\Gamma$ in a space $L_2(\rnum^n,\cnum^N)$ which satisfies
$\Gamma^2=0$, we let $\Pi=\Gamma+\Gamma^*$, and consider
perturbations of the type $\Pi_B=\Gamma+B_1\Gamma^*B_2$ where
$B_1$ has positive real part on the range of $\Gamma^*$, $B_2$ has
positive real part on the range of $\Gamma$, and
$\Gamma^*B_2B_1\Gamma^*=0$  and $\Gamma B_1B_2\Gamma =0$. In this
case there is a (non--orthogonal) Hodge decomposition of $\mH =
L_2(\rnum^n,\cnum^N)$ into closed subspaces:
$$\mH = \nul(\Pi_B)\oplus\clos{\ran(\Gamma^*_B)}
            \oplus\clos{\ran(\Gamma)}\,.$$
The quadratic estimates and functional calculus hold for
$u\in\overline{{\ran(\Pi_B)}}=\overline{{\ran(\Gamma^*_B)}}\oplus
\overline{{\ran(\Gamma)}}$.

These results have implications for spectral projections of the
Hodge--Dirac operator $d+d^*_g$ on a compact manifold $M$ with a
Riemannian metric $g$. The operator $d+d^*_g$
    is a selfadjoint operator in the Hilbert space $\mH=L_2(M,\wedge T^*M)$,
and so there is an orthogonal decomposition
$$\mH=\nul(d+d^*_g)\oplus\mH^+_g\oplus\mH^-_g$$ where
$\mH^\pm_g$ are the positive and negative eigenspaces of
$d+d^*_g$. The projections of $\mH$ onto $\mH^\pm_g$ are ${\mathbf
E}^\pm _g = \xi^\pm(d+d^*_g)$ where the   functions $\xi^\pm :
S_\mu ^o \cup\{0\}\longrightarrow \cnum$  defined by
\begin{equation*}
\xi^\pm(z)= \left\{ \begin{array}{ll}
1    & \quad \textrm{if $\pm \re z>0$}\\
0    & \quad \textrm{if $\pm \re z\leq 0$} \\
\end{array} \right.
\end{equation*}
are holomorphic on $S^o_\mu$. The subscript $g$ denotes dependence
on the metric $g$.

If the metric is perturbed to $g+h$, then the adjoint of $d$ with
respect to the perturbed metric has the form
$d^*_{g+h}=B^{-1}d^*_gB$ for an associated positive selfadjoint
multiplication operator $B$. The perturbation result
(\ref{morefun})  can be transferred to this context, thus  giving
\begin{equation}\|{\mathbf
E}^\pm_{g+h}-{\mathbf E}^\pm_g\|\lesssim\|h\|_\infty :=
\esssup_{x\in M} |h_x| \label{riemannian}
\end{equation}
    provided $\|h\|_\infty$ is not too large, where $$|h_x|=
\sup\{|h_x(v,v)| :  v\in T_xM\, ,\, g_x(v,v)=1\}.$$ What
(\ref{riemannian}) tells us is that these eigenspaces depend
continuously  on $L_\infty$ changes in the metric. Indeed the
eigenspaces depend analytically  on $L_\infty$ changes in the
metric. This  result is possibly surprising in that the local
formula for $d^*_{g+h}$ in terms of $d^*_g$ depends on the first
order derivatives of $h$.

\subsection{Acknowledgments}

This research was mostly undertaken at the Centre for Mathematics
and its Applications at the Australian National University, and
was supported by the Australian Research Council. The second
author held a visiting position at the School of Mathematics at
the University of New South Wales during the final preparation of
this paper, and thanks them for their hospitality.

   We wish to acknowledge the contributions of Pascal Auscher
and Andrea Nahmod to the development of the connections between
the Kato square root problem and quadratic estimates for a
corresponding perturbed Dirac operator. The framework developed in
their joint  paper with the third author \cite{AMcN} is a
forerunner of that presented in the current paper.

The key results of this paper were first presented at the
Conference {\it Analyse Harmonique et ses Applications} at Orsay
in June 2003, in honour of Raphy Coifman and Yves Meyer for their
profound contributions to the theory of singular integrals and to
the Calder\'on program.

\section{Statement of results}  \label{statement}

We begin by  standardizing notation and terminology. All theorems
and results in this paper are quantitative, in the sense that
constants in estimates  depend only on constants quantified in the
relevant hypotheses. Such dependence will usually be clear. We use
the notation $a \approx b$ and $b \lesssim c$, for $a,b,c\ge 0$,
to mean that there exists $C>0$ so that $ a/C \le b \le Ca$ and $b
\le C c$, respectively. The value of  $C$  varies from one usage
to the next, but then is always fixed, and depends only on
constants quantified in the relevant preceding hypotheses.

For an unbounded linear operator $A: \dom(A)\longrightarrow \mH_2$
from a domain $\dom(A)$ in a Hilbert space $\mH_1$ to another
Hilbert spaces $\mH_2$, we denote its null space by $\nul(A)$ and
its range by $\ran(A)$. The operator $A$ is said to be closed when
its graph is a closed subspace of $\mH_1\times \mH_2$. The space
of all bounded linear operators from $\mH_1$ to $\mH_2$ is denoted
$\mL(\mH_1,\mH_2)$, while $\mL(\mH):=\mL(\mH,\mH)$. See for
example \cite{Kato} for more details.

Consider three operators $\{\Gamma, B_1,B_2 \}$ in a Hilbert space
$\mH$ with the following properties.
\begin{itemize}
\item[(H1)] The operator $\Gamma:\dom(\Gamma)\longrightarrow\mH$
is a {\em nilpotent} operator from $\dom(\Gamma)\subset \mH$ to
$\mH$, by which we mean  $\Gamma$ is closed, densely defined and
$\ran(\Gamma)\subset\nul(\Gamma)$. In particular, $\Gamma^2=0$ on
$\dom(\Gamma)$. \item[(H2)] The operators $B_1$, $B_2:
\mH\longrightarrow\mH$ are bounded operators satisfying the
accretivity conditions for some $\kappa_1, \kappa_2>0$:
$$   \re (B_1 u,u) \geq \kappa_1 \|u\|^2\quad \text{for all} \quad u
\in\ran(\Gamma^*),
$$
$$  \re (B_2 u,u) \geq \kappa_2 \|u\|^2\quad \text{for all} \quad
u\in\ran(\Gamma).
$$
Let the angles of accretivity be
$$   \omega_1 := \sup_{u\in\ran(\Gamma^*)\setminus\{0\}}|\arg(B_1
u,u)| <\tfrac\pi2,$$
$$
     \omega_2 := \sup_{u\in\ran(\Gamma)\setminus\{0\}}|\arg(B_2 u
,u)| <\tfrac\pi2,$$ and set $ \omega
:=\tfrac12(\omega_1+\omega_2).$ \item[(H3)] The operators satisfy
$\Gamma^*B_2B_1\Gamma^*=0$ on $\dom(\Gamma^*)$ and $\Gamma
B_1B_2\Gamma =0$ on $\dom(\Gamma)$, that is,  $B_2B_1:
\ran(\Gamma^*) \longrightarrow \nul(\Gamma^*)$ and $B_1B_2:
\ran(\Gamma) \longrightarrow \nul(\Gamma)$. This implies that
$\Gamma B_1 ^* B_2 ^* \Gamma =0 $ on $\dom (\Gamma)$ and that
$\Gamma^* B_2 ^* B_1 ^* \Gamma^* =0$ on $\dom (\Gamma^*)$.
\end{itemize}

In some applications, $B_2$ satisfies the accretivity condition on
all of $\mH$ and $B_1={B_2}^{-1}$. In this case (H3) is
automatically satisfied, and the accretivity condition for $B_1$
holds with $\omega_1=\omega_2$.

\begin{defn} Let $\Pi = \Gamma + \Gamma^*$. Also let $\Gamma^*_B = B_1
\Gamma^* B_2 $ and $ \Gamma _B = B_2 ^* \Gamma B_1 ^*$ and then
let $\Pi_B= \Gamma + \Gamma^* _B $ and $ \Pi^*_B= \Gamma^* +
\Gamma_B.$
\end{defn}
In Section \ref{opthsec}, specifically in Lemma \ref{nilpotent}
and Corollary \ref{duality}, we show that $\Gamma^*_B =
(\Gamma_B)^*$ and $\Pi^*_B = (\Pi_B)^*$, that each of these
operators is closed and densely defined, and morever that
$\Gamma_B$ and $\Gamma^*_B$ are nilpotent.  The proofs of the
following two propositions are also given in Section
\ref{opthsec}. The first establishes a Hodge decomposition for the
perturbed operators.

\begin{prop}  \label{hodgedec}
     The Hilbert space $\mH$ has the following Hodge decomposition into
closed subspaces:
\begin{equation} \label{Hodge}
     \mH = \nul(\Pi_B)\oplus\clos{\ran(\Gamma^*_B)}
            \oplus\clos{\ran(\Gamma)} \ .
\end{equation}
Moreover, we have $\nul(\Pi_B)=\nul(\Gamma^*_B)\cap\nul(\Gamma)$
and $\clos{\ran(\Pi_B)}=\clos{\ran(\Gamma^*_B)}\oplus
\clos{\ran(\Gamma)}$. When $B_1=B_2=I$ these decompositions are
orthogonal, and in general the decompositions are topological.
Similarly, there is also a decomposition $$\mH = \nul({\Pi^*
_B})\oplus\clos{\ran(\Gamma_B)}
            \oplus\clos{\ran(\Gamma^*)} .$$
\end{prop}

\begin{defn}  \label{projections}
The bounded projections onto the subspaces in the Hodge
decomposition (\ref{Hodge}) are denoted by $\PP^0$ onto
$\nul(\Pi_B)$, $\PP^1$ onto $\clos{\ran(\Gamma^*_B)}$
   and $\PP^2$ onto $\clos{\ran(\Gamma)}$.
When $B_1=B_2=I$, these are orthogonal projections which we
denote by $\PQ^0$, $\PQ^1$ and $\PQ^2$.
\end{defn}

We  now investigate the spectrum and resolvent estimates for the
    operator $\Pi_B$.

\begin{defn} \label{psidef} Given $ 0 \leq\omega<  \mu <
\frac\pi2$, define the
     closed and open sectors and
double sectors in the complex plane by
\begin{align*}
    S_{\omega+} &:= \sett{z\in\cnum}{|\arg z|\le\omega}\cup\{0\}\,,\\
    S_{\mu+}^o &:= \sett{ z\in\cnum}{ z\ne 0, \, |\arg z|<\mu}\,,\\
    S_{\omega} &:= S_{\omega+}\cup(-S_{\omega+})\,, \\
    S_{\mu}^o &:= S_{\mu+}^o\cup(- S_{\mu+}^o)\,.
\end{align*}

Also let $\Psi(S_\mu ^o)$ denote the collection of holomorphic
functions $\psi : S_\mu ^o \longrightarrow \cnum$ such that there
exist $L,s>0$ so that
\begin{equation*} \label{eq:psi} | \psi(z)| \le L
\frac{|z|^s}{(1+|z|^{2s})}\end{equation*} for all $ z \in
S^o_\mu$. \end{defn}

\begin{prop}  \label{typeomega}
The spectrum $\sigma(\Pi_B)$ is contained in the double sector
$S_{\omega}$. Moreover the operator $\Pi_B$ satisfies resolvent
bounds
$$  \|(\I+\tau\Pi_B)^{-1} \|\lesssim
\frac{|\tau|}{\dist(\tau,S_\omega)}$$ for all $\tau \in
\cnum\setminus S_\omega$.
\end{prop}

Such an operator is  of type $S_\omega$ as defined in
\cite{ADMc,AMcNhol}. A consequence of the above proposition is
that the following operators are uniformly bounded in $t$.

\begin{defn} \label{multiop}
For $t\in\rnum\ (t\neq0)$, define the bounded operators in $\mH$:
\begin{align*}
R_t^B&:=
(\I+it\Pi_B)^{-1}\,,\\P_t^B&:=(\I+t^2{\Pi_B}^2)^{-1}=\tfrac
   12(R_t^B+R_{-t}^B)=R_t^BR_{-t}^B\quad \text{and}\\
Q_t^B&:=t\Pi_B(\I+t^2{\Pi_B}^2)^{-1}=\tfrac
     1{2i}(-R_t^B+R_{-t}^B)\end{align*}
In the unperturbed case $B_1=B_2=I$, we write $R_t$, $P_t$ and
$Q_t$ for $R_t ^B$, $P_t ^B$ and $Q_t ^B$, respectively.
\end{defn}

For an operator with the spectral properties of Proposition
\ref{typeomega}, it is useful to know whether it satisfies
quadratic estimates and whether it has a bounded holomorphic
functional calculus.
The hypotheses (H1--3) are not enough to imply quadratic
estimates. See Remark \ref{example}.  Thus we introduce further
hypotheses which allow the use of harmonic analysis.

\begin{enumerate}
\item[(H4)] The Hilbert space is $\mH=L_2(\rnum^n;\cnum^N)$, where
$n,N\in \nnum$.
   \item[(H5)]
The operators $B_1$ and $B_2$ denote multiplication by
matrix--valued functions $B_1, B_2\in
L_\infty(\rnum^n;\mL(\cnum^N))$.
   \item[(H6)] (Localisation) The nilpotent operators $\Gamma$
and $\Gamma^*$ are first order differential operators in the sense
that if $\eta:\rnum^n\longrightarrow\cnum$ is a bounded Lipschitz
function, then multiplication by $\eta$ preserves $\dom(\Gamma)$
and $\dom(\Gamma^*)$, and the commutators
$$
      \Gamma_{\nabla\eta}:=[\Gamma,\eta \I],\qquad
\Gamma^*_{\nabla\eta}:=[\Gamma^*,\eta \I]
$$
are multiplication operators such that there exists $c>0$ so that
$$|\Gamma_{\nabla\eta}(x)|, \ |\Gamma^*_{\nabla\eta}(x)| \le c
|\nabla\eta(x)|$$ for all $x\in\rnum^n$. \item[(H7)]
(Cancellation) We have $\int_{\rnum^n} \Gamma u=0$ for all
compactly supported $u\in\dom(\Gamma)$, and we have
$\int_{\rnum^n} \Gamma^* v=0$ for all compactly supported
$v\in\dom(\Gamma^*)$.
\item[(H8)] (Coercivity) There exists
$c>0$ such that
$$\|\nabla u\|\leq c \|\Pi u\|$$ for all $u\in \ran(\Pi) \cap \dom (\Pi)$.
\end{enumerate}

Observe that (H6--7) automatically hold if $\Gamma$ is a
homogeneous first order
   differential operator with constant coefficients. We now state
the first main result of the paper.

\begin{thm} \label{mainthm}
Consider the  operator $\Pi_B=\Gamma+B_1\Gamma^*B_2$ acting in the
Hilbert space $\mH=L_2(\rnum^n;\cnum^N)$, where $\{\Gamma, B_1,
B_2\}$ satisfies the hypotheses (H1--8). Then $\Pi_B$ satisfies
the quadratic estimate
\begin{equation}
\int_0^\infty\|Q_t^B u\|^2 \, \frac{dt}t=
\int_0^\infty\|\Pi_B(\I+t^2{\Pi_B}^2)^{-1}u\|^2\, t\, dt\ \approx\
\|u\|^2 \label{eq:quad}
\end{equation} for all $u \in
\clos{\ran(\Pi_B)}\subset L_2(\rnum^n;\cnum^N)$.
\end{thm}

Let us now discuss the holomorphic functional calculus for
$\Pi_B$. As a result of Proposition \ref{typeomega}, one can
define the operator $\psi(\Pi_B):\mH\longrightarrow\mH$ whenever
$\psi\in\Psi(S_\mu ^o)$ for some $\mu>\omega$, in such a way that
the mapping $\psi\mapsto\psi(\Pi_B)$ is an algebra homomorphism.
This can be done as in the Dunford functional calculus by  a
contour integral
\begin{equation}  \label{dunfordschwarz}
     \psi(\Pi_B):=\frac 1{2\pi i}\int_\gamma
\psi(\lambda)(\lambda \I-\Pi_B)^{-1} d\lambda
\end{equation}
where \label{gam} $\gamma$ is the unbounded contour $\sett{\pm r
e^{\pm i\theta}}{r\ge 0}$, $\omega <\theta<\mu$, parametrised
counterclockwise around $S_{\omega}$. The decay estimate on $\psi$
and the resolvent bounds of Proposition \ref{typeomega} guarantee
that the  integral is absolutely convergent and that $\psi(\Pi_B)$
is bounded. See for example \cite{ADMc,AMcNhol,cowl,Mc} for a
discussion of these matters.

\begin{rem}\label{otherpsi} We note in passing that each
 $\psi\in\Psi(S_\mu ^o)$ which is nonzero on both sectors
 defines a quadratic seminorm on $\mH$, and
that they are all equivalent. In particular, we have
$\int_0^\infty\| \psi(t\Pi_B)u\|^2\frac{dt}t\approx
\int_0^\infty\|Q_t^B u\|^2 \, \frac{dt}t$  for all $u\in\mH$.
Therefore, under hypotheses (H1--8),  we have $\int_0^\infty\|
\psi(t\Pi_B)u\|^2\frac{dt}t\approx \|(\I-\PP^0)u\|^2$ for all
$u\in\mH$.
\end{rem}

\begin{defn}  \label{funccalcdefn} Suppose
$\omega<\mu<\frac\pi2$. We say that $\Pi_B$ has a {\it bounded}
$S^o_\mu$ {\it holomorphic functional calculus} if
\begin{equation}  \label{bddfunccalcest}
     \|\psi(\Pi_B) \| \lesssim \|\psi\|_\infty:=\sup\{|\psi(z)|:
z\in S^o_\mu\}
\end{equation}
for all $\psi\in\Psi(S_\mu ^o)$.
\end{defn}

In this case one can define a bounded operator $f(\Pi_B)$ with
\begin{equation}  \label{bddfunccalcestinthm}
     \|f(\Pi_B) \| \lesssim \|f\|_\infty:=\sup\{|f(z)|:
z\in S^o_\mu\cup\{0\}\}
\end{equation}
for all bounded functions $f:S^o_\mu\cup\{0\}\longrightarrow\cnum$
which are holomorphic on $S^o_\mu$. The operator $f(\Pi_B)$ can be
defined by
\begin{equation}\label{evenmorefun} f(\Pi_B)u = f(0)\PP^0 u +
\lim_{n\to\infty}\psi_n(\Pi_B)u
\end{equation}
for all $u\in \mH$, where the functions $\psi_n\in\Psi(S^o_\mu)$
are uniformly bounded and tend locally uniformly to $f$ on
$S^o_\mu$; see \cite{ADMc,cowl}. The definition is independent of
the choice of the approximating sequence $(\psi_n)$. If $\Pi_B$
satisfies the quadratic estimate (\ref{eq:quad}) for all
$u\in\overline{\ran(\Pi_B)}$ then it has a bounded holomorphic
functional calculus. Thus we have the second main result of the
paper.

\begin{thm} \label{mainthm2}
Assume the hypotheses of Theorem \ref{mainthm} and let
$\omega<\mu<\frac\pi2$. Then $\Pi_B$ has a bounded $S^o_\mu$
holomorphic functional calculus in $L_2(\rnum^n;\cnum^N)$.
\end{thm}

A consequence of this theorem is that there is a decomposition of
$\mH$ into spectral subspaces. Let $\xi^\pm$ be the holomorphic
functions defined in the Introduction. Also let $\xi^0$ denote the
characteristic function of $\{0\}$ so that $\xi^0+\xi^++\xi^-=1$
on $S^o_\mu\cup\{0\}$ and $\xi^0(\Pi_B) = \PP^0$. By  Theorem
\ref{mainthm2}, the spectral projections ${\mathbf
E}_B^\pm=\xi^\pm(\Pi_B)$ are bounded, and by the functional
calculus, $\PP^0+{\mathbf E}_B^++{\mathbf E}_B^-=\I$. This leads
to part (i) of the Corollary below. Furthermore,  define the
function $\sgn$ by $\sgn(z)=z/\sqrt{z^2}$ when $z\in S^o_\mu$ and
$\sgn(0)=0$, so that $\sgn(z)=\xi_+(z)-\xi_-(z)$ and hence
$\sgn(\Pi_B)= {\mathbf E}_B^+-{\mathbf E}_B^-$. The boundedness of
this operator together with the Hodge decomposition implies part
(ii).

\begin{cor} \label{cortomainthm22}
Assume the hypotheses of Theorem \ref{mainthm}. Then
\begin{enumerate}
\item[\rm(i)] there is a (non--orthogonal) spectral decomposition
$$\mH=\nul(\Pi_B)\oplus {\mathbf E}_B^+\mH \oplus {\mathbf E}_B^-\mH
$$
 into spectral
subspaces of $\Pi_B$ corresponding to $\{0\}$,
$S_{\omega+}\setminus\{0\}$ and $S_{\omega-}\setminus\{0\}$,
respectively; and \item[\rm (ii)] we have $
\dom(\Gamma)\cap\dom(\Gamma^*_B)=\dom(\Pi_B) =
\dom(\sqrt{{\Pi_B}^2})$ with
$$\|\Gamma u\|+\|\Gamma^*_Bu\|\approx\|\Pi_Bu\| \approx
\|\sqrt{{\Pi_B}^2}u\|\,.
$$
\end{enumerate}
\end{cor}

\begin{rem} If $u\in {\mathbf E}_B^\pm\mH$, then $u(x,t) =
\exp(-t\sqrt{{\Pi_B}^2})u_0(x)$ is the solution of $\frac{\partial
u}{\partial t}\pm\Pi_Bu=0$ for $t\geq 0$ which equals $u_0$ when
$t=0$ and decays as $t \to\infty$. It is a consequence of 
Remark~\ref{otherpsi} with $\psi(z)=z\exp(-\sqrt{z^2})$ that
$\|u_0\|^2\approx\int_0^\infty\|\frac{\partial}{\partial
t}u(.,t)\|^2tdt$ for $u_0\in {\mathbf E}_B^\pm\mH$.

\end{rem}

In Section \ref{sect:cons} we use Theorem \ref{mainthm2}  and
Corollary \ref{cortomainthm22} to give a unified proof of many
results in the Calder\'on program, including the Kato square root
problem and the boundedness of the Cauchy operator on Lipschitz
curves and surfaces. We are not claiming that the approach adopted
here is always better  than the original proofs given by the
respective authors.
Nonetheless, we believe there is value in seeing that each of
these results can be easily derived from Theorem \ref{mainthm}.
Moreover, at the end of Section \ref{sect:cons} we apply Theorem
\ref{mainthm2} to Hodge--Dirac operators in Euclidean space, and
obtain  Theorem \ref{thm:dira}. This result is new.

Sections  \ref{sect:holo} and \ref{sect:riem} give further
consequences and developments of Theorems \ref{mainthm} and
\ref{mainthm2}. In Section \ref{sect:holo} we first demonstrate
that, under the hypotheses (H1--3), the resolvents of $\Pi_B$ vary
holomorphically with respect to perturbations in $B$, as do the
operators $\psi(\Pi_B)$ when $\psi\in\Psi(S^o_\mu)$. We use these
results in Theorem \ref{newthm2}, to show that, under all the
hypotheses (H1--8),  the bounded members of the functional
calculus of the perturbed Dirac operator, and quadratic functions,
depend holomorphically on perturbations in $B$. From this, we
deduce Lipschitz estimates on members of the functional calculus
of the perturbed Dirac operator $\Pi_B$, and also of the quadratic
estimates of $\Pi_B$, in terms of small perturbations in $B$. In
Section \ref{sect:riem} we prove and then apply these results to
Hodge--Dirac operators on compact Riemannian manifolds. This
enables us to establish Theorem \ref{thm:man1}, which gives
Lipschitz estimates for members of the functional calculus
(including spectral projections) of the Hodge--Dirac operator  on
compact manifolds in terms of $L_\infty$ changes in the metric. In
Appendix \ref{sect:proj}, we show that, under hypotheses (H1--3),
the Hodge projections also depend holomorphically on perturbations
in $B$, and calculate the derivatives of these projections.

We conclude this section with a brief outline of the idea behind
the proofs of Theorems \ref{mainthm} and \ref{mainthm2}. The
results in Section \ref{opthsec} just depend on hypotheses
(H1--3). We prove Propositions \ref{hodgedec} and \ref{typeomega},
and show how to reduce Theorems \ref{mainthm} and \ref{mainthm2}
to a particular quadratic estimate (\ref{thesqest}). In Section 5
we prove this estimate under all the hypotheses (H1--8). This can
be considered as a type of ``$T(b)$ argument''. In Section
\ref{princsec}, we separate out the principal part $\gamma_t$ of
the operator appearing as the integrand in the desired quadratic
estimate (\ref{thesqest}). This localization procedure relies on
Propositions \ref{hodgedec} and \ref{typeomega}, the off--diagonal
estimates established in Proposition \ref{pseudoloc}, and the
local Poincar\'e inequality together with the global coercivity
condition (H8). We estimate the principal part $\gamma_t$ of the
operator in Section \ref{Carlesonsec}. To do this we show that
$d\mu(x,t)=|\gamma_t(x)|^2\tfrac{dxdt}t$ is a Carleson measure,
and then apply Carleson's Theorem for Carleson measures. This
provides the desired result.

\section{Consequences} \label{sect:cons}

For Consequences \ref{oneDcauchy}--\ref{ex:diff} we employ the
following special case of Theorem~\ref{mainthm2}.
\begin{itemize}
\item[$\bullet$]
     Let $\cnum^N=V_1\oplus V_2$, where $V_1$ and $V_2$ are finite dimensional
complex Hilbert spaces, and form the  orthogonal direct sum
$L_2(\rnum^n;\cnum^N)= L_2(\rnum^n;V_1)\oplus L_2(\rnum^n;V_2)$.
\item[$\bullet$] Let $D$ and $D^*$ be adjoint homogeneous first
order partial differential operators with constant coefficients
$$
     D :L_2(\rnum^n;V_1)\longrightarrow L_2(\rnum^n;V_2), $$
    $$ D^* :L_2(\rnum^n;V_2)\longrightarrow L_2(\rnum^n;V_1),$$
such that there exists $c>0$ so that $$\|\nabla u\|\leq c\|Du\|
\quad \text{for all} \quad u\in\ran(D^*)\cap\dom(D)\,,$$
$$\|\nabla u\|\leq c\|D^*u\|  \quad \text{for all} \quad
u\in\ran(D)\cap\dom(D^*)\,.$$ \item[$\bullet$] The operators
$A_i:L_2(\rnum^n;V_i)\longrightarrow L_2(\rnum^n;V_i)$, $i=1,2$,
denote multiplication by functions $A_i\in
L_\infty(\rnum^n;\mL(V_i))$ which satisfy the accretivity
conditions
$$
     \re(A_1 D^*u,D^* u) \geq \kappa_1\|D^*u\|^2\quad \text{for all}
      \quad u\in\dom(D^*), $$
$$   \re(A_2 Du,Du) \geq\kappa_2\|Du\|^2\quad \text{for all} \quad
     u\in\dom(D),$$
for some $\kappa_1,\kappa_2>0$. Denote the angles of accretivity
by
\begin{align*}\omega_1&:= \sup_{u \in \dom(D^*)\setminus\nul(D^*)} |\arg(A_1 D^*u,D^*u)|
<\tfrac\pi2,\\ 
\omega_2&:= \sup_{u \in \dom(D) \setminus\nul(D)}|\arg(A_2 Du,Du)|
<\tfrac\pi2.\end{align*}
\end{itemize}
In the full space $L_2(\rnum^n;\cnum^N)=L_2(\rnum^n;V_1)\oplus
L_2(\rnum^n;V_2)$, consider the following operators:
$$
     \Gamma:=
     \begin{bmatrix}
        0 & 0 \\
        D & 0
     \end{bmatrix},
     \quad
     \Gamma^*:=
     \begin{bmatrix}
        0 & D^* \\
        0 & 0
     \end{bmatrix},
     \quad
     B_1:=
     \begin{bmatrix}
        A_1 & 0 \\
        0 & 0
     \end{bmatrix},
     \quad
     B_2:=
     \begin{bmatrix}
        0 & 0 \\
        0 & A_2
     \end{bmatrix}.
$$
With this choice of $\{\Gamma, B_1, B_2\}$, the  operator $\Pi_B$
and its square become
$$
     \Pi_B=
     \begin{bmatrix}
        0 & A_1D^*A_2 \\
        D & 0
     \end{bmatrix} \quad \text{and}
   \quad
     {\Pi_B}^2=
     \begin{bmatrix}
        A_1D^*A_2D & 0 \\
        0 & DA_1D^*A_2
     \end{bmatrix}\,.
$$
The operators $\Gamma$, $\Gamma^*$ and $\Gamma^*_B$ are clearly
nilpotent, with
\begin{align*}
     \clos{\ran(\Gamma)} & \subset L_2(\rnum^n;V_2) \subset
\nul(\Gamma)\quad\text{and}
\\
     \clos{\ran(\Gamma^*)},\, \clos{\ran(\Gamma^*_B)} & \subset L_2(\rnum^n;V_1)
     \subset \nul(\Gamma^*),\, \nul(\Gamma^*_B).
\end{align*}

\begin{thm}  \label{diffopmainthm}
     Assume that $\{D, A_1, A_2\}$ are as above, and suppose
$\omega_1+\omega_2<2\mu<\pi$. Then the operator $\Pi_B$ has a
  bounded $S^o_{\mu}$
holomorphic functional calculus in $L_2(\rnum^n;\cnum^N)$.
Moreover
\begin{itemize}
\item[{\rm (i)}] the operator $A_1D^*A_2D$ has a bounded
$S^o_{2\mu+}$ holomorphic functional calculus in
$L_2(\rnum^n;V_1)$; and
    \item[{\rm (ii)}] we have
$\dom((A_1D^*A_2D)^{1/2})=\dom(D)$ with the Kato square root
estimate
$$
     \|(A_1D^*A_2D)^{1/2}u\| \approx \|Du\|$$ for all $u\in\dom(D).$
\item[{\rm (iii)}] If furthermore $V_1=V_2=:V$, $D^*=-D$,
$A_1=A_2=:A$ and $\omega_1=\omega_2=\omega<\mu<\frac\pi2$, then
$iAD$ and $iDA$ have bounded $S^o_{\mu}$ holomorphic functional
calculi in $L_2(\rnum^n;V)$. In particular $\|\sgn(iAD)\|<\infty$
and $\|\sgn(iDA)\|<\infty$.
\end{itemize}
\end{thm}

\begin{proof}
The hypothesis of Theorem~\ref{mainthm} for this $\Pi_B$ is
satisfied, and thus by Theorem \ref{mainthm2}, $\Pi_B$ has a
bounded $S^o_{\mu}$ holomorphic functional calculus.

To prove (i), let $F:S^o_{2\mu+}\cup\{0\}\longrightarrow\cnum$ be
bounded and holomorphic on $S^o_{2\mu+}$, and write $
f(z):=F(z^2)$, $z\in S^o_\mu\cup\{0\}$. Then
$$
     f(\Pi_B)=\begin{bmatrix}
        F(A_1D^*A_2D) & 0 \\
        0 & F(DA_1D^*A_2)
     \end{bmatrix}
$$
satisfies $\|f(\Pi_B)\|\lesssim \|f\|_\infty=\|F\|_\infty$, and
thus $\|F(A_1D^*A_2D)\|\lesssim\|F\|_\infty$.

The Kato square root estimate in (ii) follows on applying
Corollary \ref{cortomainthm22} to $u\in\dom(D)$.

Now make the additional assumptions stated in (iii). That $iDA$
has a bounded $S^o_{\mu}$ holomorphic functional calculus in
$L_2(\rnum^n;V)$ can be seen as follows. Consider a bounded
function $f:S^o_{\mu}\cup\{0\}\longrightarrow\cnum$ holomorphic on
$S^o_{\mu}$.
    We find for $u\in L_2(\rnum^n;V)$, that
$$
f(\Pi_B)
     \begin{bmatrix}
        iAu \\
        u
     \end{bmatrix}
     =
f\left(
     \begin{bmatrix}
        0 & -ADA \\
        D & 0
     \end{bmatrix}
     \right)
     \begin{bmatrix}
        iAu \\
        u
     \end{bmatrix}
     =
     \begin{bmatrix}
        iA(f(iDA)u) \\
        f(iDA)u
     \end{bmatrix}.
$$
Thus $$\|f(iDA)\|\lesssim\|f(\Pi_B)\|\lesssim\|f\|_\infty\,.$$
Duality shows that $iAD=(iDA^*)^*$ also has a bounded $S^o_{\mu}$
holomorphic functional calculus in $L_2(\rnum^n;V)$, which
completes the proof of (iii). \qedend \end{proof}

Part (iii) can also be deduced from the quadratic estimates in
Theorem \ref{mainthm}, for they imply that $-ADAD$ and $-DADA$,
and hence $iAD$ and $iDA$, satisfy quadratic estimates.

We now consider several consequences of the above theorem.

\begin{ex}[The Cauchy singular integral on Lipschitz curves]
\label{oneDcauchy}
     Let $g:\rnum\longrightarrow\rnum$ be a Lipschitz function with
Lipschitz constant
$$L:=\sup_{x\ne y}\frac{|g(x)-g(y)|}{|x-y|}$$ and consider the
Lipschitz graph $\gamma:=\sett{z=x+ig(x)}{x\in\rnum}$ in $\cnum$.
The operator of differentiation with respect to $z\in\gamma$ can
be expressed in terms of the parameter $x\in\rnum$ as
$$
     D_\gamma u(x) := a Du(x) = (1+ig'(x))^{-1}u'(x) $$
where $a$ is the multiplication operator $a:v(x)\mapsto
(1+ig'(x))^{-1}v(x)$. Thus $iD_\gamma$ is of the form considered
in Theorem \ref{diffopmainthm}(iii) on making the identifications
$$
    \{n,V_1,V_2,D,D^*,A_1,A_2\}=\{1,\cnum,\cnum,\tfrac d{dx},-\tfrac d{dx},a,a\}\,.
$$

The Cauchy singular integral operator $C_\gamma$ on $\gamma$ is
then given as an operator on $L_2 (\rnum, \cnum)$ by (see
\cite{McQ,ADMc})
$$C_\gamma u(x):=\sgn(iD_\gamma)u(x) =
\frac i\pi\text{p.v.}
\int_\rnum\frac{u(y)}{(y+ig(y))-(x+ig(x))}(1+ig'(y))dy\,.$$

Using Theorem~\ref{diffopmainthm}(iii) we deduce that $iD_\gamma$
has a bounded $S^o_{\mu}$ holomorphic functional calculus in
$L_2(\rnum;\cnum)$ when $\arctan(L)<\mu<\frac\pi2$. In particular
   $\|C_\gamma\|<\infty$. The boundedness of the Cauchy integral
    $C_\gamma$ was first proved for small $L$ by
Calder\'on~\cite{Ca}, and in the general case by Coifman--\Mcc
Intosh--Meyer~\cite{CMcM}. Boundedness of other operators in the
functional calculus of $iD_\gamma$ have been proved by
Coifman--Meyer~\cite{CM}, Kenig--Meyer~\cite{KM} and \Mcc
Intosh--Qian~\cite{McQ}.
\end{ex}

\begin{ex}[The one dimensional Kato square root problem]
     Let $a\in L_\infty(\rnum;\cnum)$ be such that $\re a(x)\ge \kappa>0$
     for almost every $x$,
and denote the angle of accretivity by $\omega:=\esssup|\arg
a(x)|$. In Theorem~\ref{diffopmainthm}, let
$$
     \{n,V_1,V_2,D,D^*,A_1,A_2\}=\{1,\cnum,\cnum,\tfrac d{dx},-\tfrac
d{dx},I,a\}
$$
where $a$ is the multiplication operator $a:f(x)\mapsto a(x)f(x)$,
and suppose $\omega<\mu<\frac\pi2$. By
Theorem~\ref{diffopmainthm}(i) we deduce that $-\tfrac
d{dx}a\tfrac d{dx}$ has a bounded $S^o_{2\mu+}$ holomorphic
functional calculus in $L_2(\rnum;\cnum)$. This can be proved by
abstract methods since $-\tfrac d{dx}a\tfrac d{dx}$ is a maximal
accretive operator, see \cite{ADMc}. However,
Theorem~\ref{diffopmainthm}(ii) proves the Kato square root
estimate in one dimension:
\begin{equation}\label{onedimkato}
     \big\|(-\tfrac d{dx}a\tfrac d{dx})^{1/2}u\big\| \approx \big\|\tfrac
{du}{dx}\big\|\end{equation}  for all $ u\in H^1(\rnum).$ This
estimate was first proved by Coifman--\Mcc
Intosh--Meyer~\cite{CMcM}.
\end{ex}

\begin{rem} \label{example}
It is known that (\ref{onedimkato}) may fail if $D$ and $A_2$ are
not differentiation and multiplication operators \cite{Alan1}.
Working backwards, we find that hypotheses (H1--3) are not
sufficient to ensure that $\Pi_B$ satisfies quadratic estimates or
that it has a  bounded holomorphic functional calculus.
\end{rem}

\begin{ex}
Let $a_i\in L_\infty(\rnum;\cnum)$, for $i=1,2$, be such that
there exists $\kappa>0$ so that $\re a_i(x)\ge \kappa>0$ for
almost every $x$, and denote the angles of accretivity by
$\omega_i:=\esssup|\arg a_i(x)|$. In Theorem~\ref{diffopmainthm},
let
$$
     \{n,V_1,V_2,D,D^*,A_1,A_2\}=\{1,\cnum,\cnum,\tfrac d{dx},-\tfrac
d{dx},a_1,a_2\}
$$
where $a_i$ is the multiplication operator $a_i:f(x)\mapsto
a_i(x)f(x)$, and suppose $\omega_1+\omega_2<2\mu<\pi$. By
Theorem~\ref{diffopmainthm}(i) we deduce that $-a_1\tfrac
d{dx}a_2\tfrac d{dx}$ has a bounded $S^o_{2\mu+}$ holomorphic
functional calculus in $L_2(\rnum;\cnum)$. This result was first
proved by Auscher--\Mcc Intosh--Nahmod~\cite{AMcN} (though with
$\mu>\max\{\omega_1,\omega_2\}$). Further
Theorem~\ref{diffopmainthm}(ii) proves the estimate
$$
     \big\|(-a_1\tfrac d{dx}a_2\tfrac d{dx})^{1/2}u\big\| \approx \big\|\tfrac
{du}{dx}\big\|$$ for all $ u\in H^1(\rnum). $ This estimate was
first proved by Kenig--Meyer~\cite{KM}. A proof is also given in
\cite{AMcN}, using a framework which can be considered a
forerunner of the approach developed here.
\end{ex}

\begin{ex}[The Clifford--Cauchy singular integral on a Lipschitz surface]
    Let $g:\rnum^n\longrightarrow\rnum$ be a Lipschitz function with
Lipschitz constant $L$, and consider the Lipschitz graph
$\Sigma:=\sett{(x,g(x))}{x\in\rnum^n}$ in $\rnum^{n+1}$. On
identifying $\rnum^{n+1}$ with $\Lambda^0\oplus\Lambda^1$ in the
complex Clifford algebra $\cnum_{(n)} ( \approx\wedge_\cnum
\rnum^n)$ generated by $\rnum^n$, where the generating basis
$\{e_i\}$ satisfies the canonical commutation relation
$e_ie_j+e_je_i=-2\delta_{ij}$, then
$\Sigma=\sett{g(x)+x}{x\in\rnum^n}$. Furthermore, let $\dirac$
denote the Dirac operator
$$
    \dirac u(x):=\sum_{k=1}^n e_k\frac{\partial u}{\partial x_k}(x),
    \qquad u:\rnum^n\longrightarrow\cnum_{(n)}.
$$
This first order partial differential operator $\dirac$ is
elliptic and selfadjoint. In Theorem~\ref{diffopmainthm}, let
$$
\{n,V_1,V_2,D,D^*,A_1,A_2\}=\{n,\cnum_{(n)},\cnum_{(n)},-i\dirac,i\dirac,A,A\}
$$
where $A$ is the multiplication operator $A:u(x)\mapsto (1-\dirac
g(x))^{-1}u(x)$. In this case, we define the operator
$\dirac_\Sigma$ on $L_2(\rnum^n,\cnum_{(n)})$ by
$$   \dirac_\Sigma u(x) := A\dirac u(x) = (1-\dirac g(x))^{-1}\dirac
u(x)
$$
and, parametrizing $\Sigma$ with $g(x)+x$, the Cauchy singular
integral operator $C_\Sigma$ on $\Sigma$ is given by
\begin{align*}
      C_\Sigma u(x) & :=\sgn(\dirac_\Sigma)u(x) \\ &\ =
     \frac
2{\sigma_n}
\text{p.v.}\int_{\rnum^n}\frac{(g(x)-x)-(g(y)-y)}{(|y-x|^2+(g(y)-g(x))^2)^
{(n+1)/2}} (1-\dirac g(y))u(y)dy
\end{align*}
where $\sigma_n$ is the volume of the unit $n$-sphere in
$\rnum^{n+1}$.

Suppose $\omega:=\arctan(L)<\mu<\frac\pi2$. By
Theorem~\ref{diffopmainthm}(iii) we deduce that $ \dirac_\Sigma$
has a bounded $S^o_{\mu}$ holomorphic functional calculus in
$L_2(\rnum^n;\cnum_{(n)})$, and in particular that
$\|C_\Sigma\|<\infty$. The boundedness of the Clifford--Cauchy
integral $C_\Sigma$ follows from the boundedness of the Cauchy
integral in Consequence~\ref{oneDcauchy} using Calder\'on's
rotation method (c.f. \cite{CMcM}). A direct proof of the
boundedness of $C_\Sigma$ using Clifford analysis was first given
by Murray~\cite{Mu} for surfaces with small $L$, and in the
general case by \Mcc Intosh~\cite{Mc1}. Boundedness of  the
functional calculus of $\dirac_\Sigma$ has been proved by Li--\Mcc
Intosh--Semmes~\cite{LMcS} and Li--\Mcc
Intosh--Qian~\cite{LMcQ}.
\end{ex}

In the following three consequences,  the differential operator
$D$ no longer has dense range.

\begin{ex}[The Kato square root problem] \label{ex:kato}
     Let $A\in L_\infty(\rnum^n;\mL(\cnum^n))$ be such that
$\re (A(x)v,v)\ge \kappa>0$ for every $v\in\cnum^n$ with $|v|=1$,
and almost every $x$, and denote the angle of accretivity by
$\omega:=\esssup_{v,x}|\arg (A(x)v,v)|$. In
Theorem~\ref{diffopmainthm}, let
$$
     \{n,V_1,V_2,D,D^*,A_1,A_2\}=\{n,\cnum,\cnum^n,\nabla,-\text{div},I,A\}
$$
where $A$ denotes the multiplication operator $A:u\mapsto Au$,
and suppose $\omega<\mu<\frac\pi2$.
   From Theorem~\ref{diffopmainthm}(i) we deduce that
$-\text{div}A\nabla$ has a bounded $S^o_{\mu+}$ holomorphic
functional calculus in $L_2(\rnum^n;\cnum)$. This can be proved by
abstract methods since $-\text{div}A\nabla$ is a maximal accretive
operator, see \cite{ADMc}. More importantly,
Theorem~\ref{diffopmainthm}(ii) implies the full Kato square root
estimate
$$
     \big\|(-\text{div}A\nabla)^{1/2}u\big\| \approx \big\|\nabla
u\big\|$$ for all $ u\in H^1(\rnum^n). $ This result was
proved in a series of papers by Hofmann--\Mcc Intosh~\cite{HMc}, Auscher--Hofmann--Lewis--Tchamitchian~\cite{AHLT}, Hofmann--Lacey--\Mcc Intosh \cite{HLMc}, and, in full generality, by Auscher--Hofmann--Lacey--\Mcc
Intosh--Tchamitchian~\cite{AHLMcT}. Earlier
results on the Kato square root problem are due to
Fabes--Jerison--Kenig~\cite{FJK} and
Coifman--Deng--Meyer~\cite{CDM}, where $A$ is assumed to be close
to the identity, and to \Mcc Intosh~\cite{Mc2} when H\"older
continuity of $A$ is assumed. For many more partial results, see
the book of Auscher and Tchamitchian \cite{ATch}. This book
provides an important bridge between the one--dimensional results
and the current theory.
\end{ex}

\begin{ex}\label{ex:nahmod}
Let $a\in L_\infty(\rnum^n;\cnum)$ be such that $\re a(x)\ge
\kappa>0$ for almost every $x$, and let $A\in
L_\infty(\rnum^n;\mL(\cnum^n))$ be such that $\re (A(x)v,v)\ge
\kappa>0$ for every $v\in\cnum^n$, $|v|=1$, and almost every $x$.
Denote the angles of accretivity by $\omega_1:=\esssup|\arg a(x)|$
and $\omega_2:=\esssup_{v,x}|\arg (A(x)v,v)|$. In
Theorem~\ref{diffopmainthm}, let
$$
     \{n,V_1,V_2,D,D^*,A_1,A_2\}=\{n,\cnum,\cnum^n,\nabla,-\text{div},a,A\}
$$
where $a$ is the multiplication operator $a:u(x)\mapsto a(x)u(x)$
and $A$ is the multiplication operator $A:v(x)\mapsto A(x)v(x)$.
   From Theorem~\ref{diffopmainthm}(i) we deduce that
$-a\,\text{div}A\nabla$ has a bounded $S^o_{2\mu+}$
     holomorphic functional calculus in
$L_2(\rnum^n;\cnum)$ when $\omega_1+\omega_2<2\mu<\pi$. This was
proved by \Mcc Intosh--Nahmod~\cite{McN} in the case when $A=I$,
and by  Duong--Ouhabaz~\cite{DO} under regularity assumptions on
$A$. Theorem~\ref{diffopmainthm}(ii) also proves the estimate
$$
     \big\|(-a\,\text{div}A\nabla)^{1/2}u\big\| \approx \big\|\nabla
u\big\|$$ for all $ u\in H^1(\rnum^n). $
\end{ex}

\begin{ex}[The Kato square root problem for systems]
\label{katosys} Let $W$ be a f\-in\-i\-t\-e dimensional Hilbert
space and let $A\in L_\infty(\rnum^n;\mL(\cnum^n\otimes W))$ be
such that
$$   \re \int_{\rnum^n}(A(x)\nabla u(x),\nabla u(x))\,dx\ge \kappa\|\nabla
u\|^2$$ for all $u\in H^1(\rnum^n;W)$ and some $\kappa>0$.
In Theorem~\ref{diffopmainthm}, let
$$
     \{n,V_1,V_2,D,D^*,A_1,A_2\}=\{n,W,\cnum^n\otimes W,\nabla,-\text{div},I,A\}
$$
where $A$ is the multiplication operator $A:f(x)\mapsto A(x)f(x)$.
Theorem~\ref{diffopmainthm}(ii) proves the Kato square root
estimate for these elliptic systems:
$$
     \big\|(-\text{div}A\nabla)^{1/2}u\big\| \approx \big\|\nabla
u\big\|$$ for all $u\in H^1(\rnum^n;W). $ This estimate was first
proved by Auscher--Hofmann--\Mcc
Into\-s\-h--Tchamitchian~\cite{AHMcT}.
\end{ex}

\begin{ex}[Differential forms] \label{ex:diff}
For $n \geq1$, let $\Lambda = \oplus_{i=0}^n
\Lambda^i=\wedge_\cnum\rnum^n$ denote the complex exterior algebra
over $\rnum^n$. Let $B$ be a bounded multiplication operator on
$L_2(\rnum^n; \Lambda )$ with bounded inverse which satisfies the
following accretivity condition: there exists $\kappa>0$ such that
for almost every $x \in \rnum^n$, we have
$$ \re(B(x) v,v) \ge
\kappa | v |^2$$ for every $v \in \Lambda$. Let $d$ denote the
exterior derivative,  and consider the perturbed Hodge--Dirac
operator $D_B = d + B^{-1} d^* B$. We further suppose that $B$
splits over $L_2(\rnum^n, \Lambda^0) \oplus \dots \oplus
L_2(\rnum^n, \Lambda^n)$ as $B^0 \oplus \dots \oplus B^n$, and so
$D_B$ can be illustrated by the following diagram.

$$\begin{array}{ccccccc}
    L_2(\rnum^n,\Lambda^0)& \stackrel{d=\nabla}{\bf
\longrightarrow} &L_2(\rnum^n,\Lambda^1)& \stackrel{d}{\bf
\longrightarrow}& \dots & \stackrel{d}{\longrightarrow}
&L_2(\rnum^n,\Lambda^n)
\\[3mm]
{\bf \Big\downarrow}B^0&&{\bf \Big\downarrow}
B^1&& &&  \Big\downarrow B^n\\[5mm]
    L_2(\rnum^n,\Lambda^0)&
\stackrel{d^*=-\divv}{\longleftarrow} & L_2(\rnum^n,\Lambda^1)&
\stackrel{d^*}{\longleftarrow}& \dots &
\stackrel{d^*}{\longleftarrow} & L_2(\rnum^n,\Lambda^n)
\end{array}$$

Let $\omega>0$ denote the angle of accretivity of $B$, and let $
\omega < \mu < \frac\pi2$. We now apply Theorem \ref{mainthm2} and
Corollary \ref{cortomainthm22} with $\Gamma =d$, $B_1:=B^{-1}$ and
$B_2:=B$, to obtain the following new result.

\begin{thm} \label{thm:dira}
The operator $D_B$  has a bounded $S^o_\mu$ holomorphic functional
calculus in $L_2(\rnum^n, \Lambda)$. Moreover, the operator
   ${D_B}^2$ has a bounded $S^o_{2\mu+}$ holomorphic
functional calculus in $L_2(\rnum^n, \Lambda )$. Furthermore,
$\dom(d) \cap \dom(d^* B)= \dom(\sqrt{{D_B}^2 })$ with
$$ \|d u \| + \| d^* B u \| \approx \left\| \sqrt{ {D_B}^2 }u \right\|\,.$$
\end{thm}

The restriction of the second and third claims to $u\in
L_2(\rnum^n,\Lambda^0)$ provides an alternative approach to the
results obtained in Consequences \ref{ex:kato} and
\ref{ex:nahmod}, though not those of Consequence \ref{katosys}.
The implications for the full exterior algebra are new and will be
developed further in Remark \ref{deeper}.
\end{ex}

\section{Operator theory of $\Pi_B$}  \label{opthsec}

Throughout this section, we assume that the triple of operators
$\{\Gamma, B_1,B_2 \}$ in a Hilbert space $\mH$ satisfies
properties (H1--3). We prove Propositions~\ref{hodgedec} and
\ref{typeomega}, and then show how to reduce
Theorems~\ref{mainthm} and \ref{mainthm2} to a quadratic estimate
which will be proved in Section~\ref{harmansec}.

Let us start by recording the following useful consequences of
(H2):
\begin{align}
\|B_1u\|&\approx\|u\|\approx\|B_1^*u\|\ \quad \text{for all}\quad
u\in
\overline{\ran(\Gamma^*)}\,;\label{equiv1}\\
\|B_2u\|&\approx\|u\|\approx\|B_2^*u\|\ \quad \text{for all}\quad
u\in \overline{\ran(\Gamma)}\,.\label{equiv2}
\end{align}

\begin{lem} \label{nilpotent} The operators $\Gamma_B^*:=B_1\Gamma^*B_2$ and
$\Gamma_B:=B_2^*\Gamma B_1^*$ are nilpotent, and
$(\Gamma_B)^*=\Gamma_B^*$.
\end{lem}

\begin{proof}
 First note that by (H3), $\ran(\Gamma^*_B)\subset \nul(\Gamma^*_B)$
 and $\ran(\Gamma_B)\subset \nul(\Gamma_B)$.
 To prove that the two operators are densely defined,
 closed and adjoint, we use the following operator theoretic
 fact: Let $A$ be a closed and densely defined operator
 and let $T$ be a bounded operator.
 Then $TA$ is densely defined, $A^*T^*$ is closed
 and $(TA)^*=A^*T^*$. If furthermore $\|Tu\|\approx\|u\|$ for
 all $u\in\ran (A)$, then $TA$ and $A^*T^*$ are closed,
 densely defined and adjoint operators.
 Applying this fact first with $A=\Gamma^*$, $T=B_1$ and
 then with $A=\Gamma B_1^*$, $T=B_2^*$ proves the lemma.
\qedend \end{proof}

We next prove a lemma concerning the operators
$\Pi_B:=\Gamma+\Gamma^*_B$ with
$\dom(\Pi_B)=\dom(\Gamma)\cap\dom(\Gamma^*_B)$, and
$\Pi_B^*:=\Gamma^*+\Gamma_B$ with
$\dom(\Pi^*_B)=\dom(\Gamma^*)\cap\dom(\Gamma_B)$.

\begin{lem} \label{lem:gam1}  We have
\begin{align*}\| \Gamma u \| +\| \Gamma_B^* u \|  &\approx \,  \| \Pi_B u
\| \quad \text{for all} \quad u \in \dom(\Pi_B), \quad
\text{and}\\
\| \Gamma^* u \| + \| \Gamma_B u \| &\approx \,  \| \Pi_B^* u \|
\quad \text{for all} \quad u \in \dom(\Pi_B^*).\end{align*}
\end{lem}

\begin{proof}
The first estimate follows from the observation that (H2--3)
implies
$$ \| \Gamma
u \|^2 \lesssim |( B_2 \Gamma u , \Gamma u )| = | (B_2 \Pi_B
u,\Gamma u )| \lesssim \| \Pi_B u \| \, \| \Gamma u \|$$ for every
$u \in \dom(\Pi_B)$. The other claims follow by similar reasoning.
\qedend \end{proof}

We now prove Proposition \ref{hodgedec} and then Proposition
\ref{typeomega}.

\begin{proof}[Proof of Proposition~\ref{hodgedec}]
It is an immediate consequence of the lemma that
$\nul(\Pi_B)=\nul(\Gamma^*_B)\cap\nul(\Gamma)$.

     Note that once we prove
\begin{equation} \label{prehodge}
     \mH=\clos{\ran(\Gamma^*_B)}\oplus\nul(\Gamma)=\nul(\Gamma^*_B)\oplus
\clos{\ran(\Gamma)}
\end{equation}
then the Hodge decomposition follows since
$\clos{\ran(\Gamma^*_B)}\subset \nul(\Gamma^*_B)$ and
$\clos{\ran(\Gamma)}\subset\nul(\Gamma)$ by nilpotence. In the
case $B_1=B_2=I$, (\ref{prehodge}) is orthogonal since $\Gamma$
and $\Gamma^*$ are adjoint operators. To prove (\ref{prehodge})
for a general $B$, it suffices to prove the four statements
$$\mH\supset\clos{\ran(\Gamma^*_B)}\oplus\nul(\Gamma),  \quad
\mH\supset\nul(\Gamma^*_B)\oplus \clos{\ran(\Gamma)}, $$ $$
\mH\supset\clos{\ran(\Gamma^*)}\oplus\nul(\Gamma_B), \quad
\mH\supset\nul(\Gamma^*)\oplus \clos{\ran(\Gamma_B)},$$ and use
duality.

Let us  consider the first of these. We need to show that
$$\|\Gamma^*_Bu\|+\|v\|\lesssim \|\Gamma^*_Bu + v\|$$ for all
$u\in\dom(\Gamma^*_B)=\dom(\Gamma^*B_2)$ and $v\in\nul(\Gamma)$.
This follows from
\begin{align*}
 \| \Gamma^*B_2 u\|^2  \lesssim |\re(B_1\Gamma^*B_2 u,
\Gamma^*B_2u)| &= |\re(\Gamma^*_B u+v, \Gamma^*B_2u)| \\ &\leq
\|\Gamma^*_B u +v\| \| \Gamma^*B_2 u\|.\end{align*}

For the second statement we need to show that $$\| v \|+\| \Gamma
u\|\lesssim \| v + \Gamma u\|$$ for all $u\in\dom(\Gamma)$ and
$v\in \nul(\Gamma^*_B)=\nul(\Gamma^*B_2)$. This follows from
$$\|\Gamma u\|^2 \lesssim|(\Gamma u,B^*_2\Gamma u)| = |(v+\Gamma
u,B^*_2\Gamma u)|\lesssim \|v+\Gamma u\| \|\Gamma u \|.$$

The third and fourth statements have  similar proofs. \qedend
\end{proof}

\begin{cor} \label{duality} The operators $\Pi_B$ and $\Pi^*_B$ are
closed, have dense  domains, and satisfy $(\Pi_B)^*=\Pi^*_B$.
\end{cor}

This is a straightforward consequence of the preceding results. We
are now in a position to prove the spectral properties stated in
Section 2.

\begin{proof}[Proof of Proposition~\ref{typeomega}]
Let $f=(\I+\tau\Pi_B)u$ where $\tau\in\cnum\setminus S_\omega$ and
$u\in\dom(\Pi_B)$. To prove the estimate $\|u\|\lesssim\|f\|$, use
Proposition~\ref{hodgedec} to write
$$
     f=f_0+f_1+f_2,\, u=u_0+u_1+u_2\in
     \nul(\Pi_B)\oplus\clos{\ran(\Gamma^*_B)}
            \oplus\clos{\ran(\Gamma)}
$$
and $f_1=B_1\tilde f_1$, $u_1=B_1\tilde u_1$, where $\tilde f_1$,
$\tilde u_1\in \clos{\ran(\Gamma^*)}$. We obtain the system of
equations
\begin{align*}   f_0 &= u_0\\
f_1 &= u_1+ \tau\Gamma^*_B u_2, \text{ thus by (\ref{equiv1}),}\
    \tilde f_1 = \tilde
u_1+ \tau\Gamma^* B_2 u_2\\
     f_2 &= u_2+ \tau\Gamma u_1\,.
\end{align*}
These equations imply the identity
\begin{equation}\label{identity}
     -\conj \tau(\tilde u_1,B_1\tilde u_1)+\tau(B_2u_2,u_2)= -\conj \tau(\tilde
f_1,B_1 \tilde u_1)+\tau(B_2u_2,f_2)\, .
\end{equation}
Let $$\theta_1=\arg(\tilde u_1,B_1\tilde u_1),\quad \text{and}
\quad \theta_2=\arg(B_2u_2,u_2)$$ so that by (H2),
$|\frac12\theta_1-\frac12\theta_2|\leq\omega$. Suppose for a
moment that $\im \tau>0$ and let $\mu=\arg\tau$. Then
\begin{align}
| -\conj \tau&(\tilde u_1,B_1\tilde u_1)+\tau(B_2u_2,u_2)|\nonumber\\
&\geq \im e^{-i(\theta_1+\theta_2)/2}\left(-\conj\tau(\tilde
u_1,B_1\tilde u_1)+\tau(B_2u_2,u_2)
\right)\nonumber\\
&=|\tau| \sin(-\tfrac12\theta_1+\tfrac12\theta_2
+\mu)\left(|(\tilde u_1,B_1\tilde u_1)|+ |(B_2 u_2,
u_2)|\right)\label{ineq}\\
&\geq \dist(\tau,S_\omega)\left( |(\tilde u_1,B_1\tilde
u_1)|+|(B_2 u_2, u_2)|\right)\ .\nonumber
\end{align}
Therefore, by (H2), (\ref{identity}) and (\ref{ineq}),
$$\|\tilde u_1\|^2+ \|u_2\|^2\lesssim
|(\tilde u_1,B_1\tilde u_1)|+|(B_2 u_2, u_2)| \lesssim
\frac{|\tau|}{\dist(\tau,S_\omega)}(\|\tilde f_1\|\|\tilde u_1\|+
\|u_2\|\|f_2\|)$$ and thus
$$\|u\|\approx\|u_0\|+\|u_1\|+\|u_2\|\lesssim\frac{|\tau|}{\dist
(\tau,S_\omega)}\|f\|.$$ A slight variation gives the estimate for
$\im\tau<0$.

Finally, applying the proof above to
$\I+\overline\tau\Pi^*_B=(\I+\tau\Pi_B)^*$ shows that
$\I+\tau\Pi_B$ is surjective. \qedend \end{proof}

\begin{cor} The operator ${\Pi_B}^2 =\Gamma B_1\Gamma^*B_2 +
B_1\Gamma^*B_2\Gamma$ is closed, has dense domain, its spectrum
$\sigma({\Pi_B}^2)$ is contained in the sector $S_{2\omega+}$, and
it satisfies resolvent bounds
$\|(\I-\tau^2{\Pi_B}^2)^{-1}\|\lesssim
\frac{|\tau^2|}{\dist(\tau^2,S_{2\omega+})}$ for all $\tau \in
\cnum\setminus S_\omega$.
\end{cor}

Such an operator is said to be of type $S_{2\omega+}$ in
\cite{ADMc} and of type $2\omega$ in \cite{AMcNhol,cowl}.

\begin{rem} \label{commutes} Note that $\Pi_B$ intertwines $\Gamma$ and
$\Gamma^*_B$ in the sense that $\Pi_B\Gamma u=\Gamma^*_B\Pi_B u$
for all $u\in \dom(\Gamma^*_B\Pi_B)$ and $\Pi_B\Gamma^*
u=\Gamma\Pi_B u$ for all $u\in \dom(\Gamma\Pi_B)$. Thus
${\Pi_B}^2$ commutes with both $\Gamma$ and $\Gamma^*_B$ on the
appropriate domains. We find that $\Gamma P^B_tu= P^B_t\Gamma u$
for all $u \in\dom(\Gamma)$ and  $\Gamma^*_B P^B_tu=
P^B_t\Gamma^*_B u$ for all $u \in\dom(\Gamma^*_B)$.
\end{rem}

We saw in Definition \ref{multiop} that the operators $P^B_t$ and
$Q^B_t$ are uniformly bounded in $t$. A consequence of this is the
identity
\begin{equation}  \label{qtresolution}
     \int_0^\infty (Q^B_t)^2u\,\frac{dt}t =
\lim_{\substack{\alpha\to0\\ \beta\to\infty}} \int_\alpha^\beta
(Q^B_t)^2u\,\frac{dt}t =\tfrac12 \lim_{\substack{\alpha\to0\\
\beta\to\infty}}(P^B_\alpha-P^B_\beta)u=\tfrac12 (\I-\PP^0) u
\end{equation}
for all $u\in\mH$. (Verify this on $\nul(\Pi_B)$ and for
$u\in\dom(\Pi_B)\cap\ran(\Pi_B)$ which is dense in
$\overline{\ran(\Pi_B)}$ and use the uniform boundedness.)
   For the selfadjoint operator
$\Pi$ this can be proved by the usual spectral theory, and has the
following consequence.
\begin{lem}  \label{schur}
The  quadratic estimate
\begin{equation}  \label{standardsqest}
     \int_0^\infty\|Q_t u\|^2\, \frac{dt}t\leq\tfrac12\|u\|^2
\end{equation}
holds for all $u\in\mH$. \end{lem}

We use the following operator in the proof of Theorem
\ref{mainthm}.

\begin{defn} \label{defnope}
Define, for all  $t\in\rnum$, the bounded operators
$$\Theta^B _t:=t\Gamma^*_B(\I+t^2 {\Pi_B}^2
)^{-1}\,.$$
\end{defn}
By Remark \ref{commutes},  $\Theta^B _tu= (\I+t^2 {\Pi_B}^2
)^{-1}t\Gamma^*_Bu$ for all $u\in\dom(\Gamma^*_B)$, and
consequently $\Theta^B _tu=Q^B_tu$ for all $u\in\nul(\Gamma)$.

\begin{prop}\label{whatisneeded}
Consider the  operator $\Pi_B=\Gamma+B_1\Gamma^*B_2$ acting in a
Hilbert space $\mH$, where $\{\Gamma, B_1, B_2\}$ satisfies the
hypotheses (H1--3). Also assume that the estimate
\begin{equation} \label{thesqest}
     \int_0^\infty\|\Theta^B _t P_t u\|^2\,\frac{dt}t\leq c\|u\|^2
\end{equation}
   holds for all $u\in\ran(\Gamma)$ and some constant $c$,
together with the three similar estimates obtained on replacing
$\{\Gamma,B_1, B_2\}$ by $\{\Gamma^* ,B_2, B_1\}$,
$\{\Gamma^*,B_2^*, B_1^*\}$ and $\{\Gamma,B_1^*, B_2^*\}$.
   Then $\Pi_B$ satisfies the
quadratic estimate (\ref{eq:quad}) for all
$u\in\overline{\ran(\Pi_B)}$, and has a bounded holomorphic
$S^o_\mu$ functional calculus.
\end{prop}

\begin{proof}\label{proof:thm}  (i) We start by proving the estimate
\begin{equation}  \label{firsttermprop}
      \int_0^\infty\|\Theta^B _t(\I-P_t) u\|^2\,\frac{dt}t\lesssim \|u\|^2
\end{equation}
for all $u\in \ran(\Gamma)$. We use the orthogonal projection
$\PQ^2:\mH\longrightarrow\clos{\ran(\Gamma)}$ and the bounded
projection $\PP^1:\mH\longrightarrow\clos{\ran(\Gamma^*_B)}$.
Since $u\in\ran(\Gamma)$ implies $P_t u\in\ran(\Gamma)$, we obtain
$$
      \Theta^B _t(\I-P_t) u = \Theta^B _t \PQ^2 (\I-P_t) u
      =Q^B _t t \Gamma Q_t  u= (\I-P^B _t) \PP^1 Q_t u
$$
and thus $\|\Theta^B _t(\I-P_t) u\| \lesssim \| Q_t u\|$ for all
$u\in\ran(\Gamma)$. This with (\ref{standardsqest}) proves
(\ref{firsttermprop}).

We remark that this use of the Hodge decompositions to handle the
$(\I-P_t)$ term is a key step in the  proof of Theorem
\ref{mainthm}.

(ii) A combination of (\ref{thesqest}) with (\ref{firsttermprop})
gives
   the estimate
\begin{equation} \label{gettingthere}
\int_0^\infty\|Q^B _t u\|^2\,\frac{dt}t= \int_0^\infty\|\Theta^B
_t u\|^2\,\frac{dt}t\lesssim\|u\|^2
\end{equation}
for all $u\in \overline{\ran(\Gamma)}$.

Now the hypotheses of the theorem remain unchanged  on replacing
$\{\Gamma,B_1, B_2\}$ by $\{\Gamma^* ,B_2, B_1\}$, in which case
the estimate in (\ref{gettingthere}) becomes
$$
     \int_0^\infty\|tB_2\Gamma B_1(\I+t^2(\Gamma^*+B_2\Gamma B_1)^2)^{-1}
v\|^2\,\frac{dt}t\lesssim\|v\|^2
$$ for all $v\in\clos{\ran(\Gamma^*)}
$. Using the assumption $\Gamma B_1B_2\Gamma =0$, we get
$$\Gamma B_1(\I+t^2(\Gamma^*+B_2\Gamma B_1)^2)^{-1}=\Gamma
(\I+t^2 \Pi_B ^2 )^{-1}B_1$$ and thus, by (\ref{equiv1}) and
(\ref{equiv2}),
\begin{align*}
     \int_0^\infty\|t\Gamma(\I+t^2 \Pi_B ^2 )^{-1} B_1v\|^2\,\frac{dt}t&\lesssim
     \int_0^\infty\|tB_2\Gamma B_1(\I+t^2(\Gamma^*+ B_2\Gamma B_1) ^2
)^{-1} v\|^2\,\frac{dt}t\\
     &\lesssim\|v\|^2\lesssim\|B_1v\|^2
\end{align*}
for all $v\in\clos{\ran(\Gamma^*)}.$ Hence
$$\int_0^\infty\|Q^B _t u\|^2\,\frac{dt}t=
\int_0^\infty\|t\Gamma(\I+t^2{\Pi_B}^2)^{-1}
u\|^2\,\frac{dt}t\lesssim\|u\|^2$$ for all $u\in
\overline{\ran(\Gamma^*_B)}$.

On recalling the Hodge decompositon $\mH =
\nul(\Pi_B)\oplus\clos{\ran(\Gamma^*_B)}
            \oplus\clos{\ran(\Gamma)}$, and noting that $Q^B_t=0$ on
$\nul(\Pi_B)$, we conclude that
$$\int_0^\infty\|Q^B _t u\|^2\,\frac{dt}t\lesssim\|u\|^2$$
for all $u\in \mH$.

(iii) To prove the reverse square function estimate, consider the
adjoint operator $\Pi_B^*=\Gamma^*+B_2^*\Gamma B_1^*$. From (ii)
applied to $\Pi_B^*$, we get
$$
     \int_0^\infty\|(Q^B_t)^* v\|^2\,\frac{dt}t \lesssim \|v\|^2$$ for all $v\in
\mH$. By (\ref{qtresolution}), we
   have the resolution of the identity
$$
     \int_0^\infty (Q^B_t)^2u\,\frac{dt}t = \tfrac{1}{2} u$$ for all $
u\in\clos{\ran(\Pi_B)}$, and thus
\begin{align*}
     \|u\|\lesssim \sup_{\|v\|=1} |(u,v)|
     & \approx
     \sup_{\|v\|=1} \left|\left(\int_0^\infty
(Q^B_t)^2u\,\frac{dt}t ,v\right)\right| \\
     & =
     \sup_{\|v\|=1} \left|\int_0^\infty (Q^B_tu
,(Q^B_t)^*v)\frac{dt}t\right| \\
     & \lesssim\left( \int_0^\infty\|Q^B_t u\|^2\,\frac{dt}t \right)^{1/2}
\end{align*}
for all $ u\in\clos{\ran(\Pi_B)}.$ This completes the proof that
(\ref{eq:quad}) holds for all $u\in\overline{\ran(\Pi_B)}$. This
procedure is standard, at least when $\nul(\Pi_B)=0$. (See e.g.
\cite{ADMc}.)

(iv) It is also well-known that quadratic estimates imply the
boundedness of the functional calculus. We include a proof for
completeness.

Note that a direct norm estimate using (\ref{dunfordschwarz})
shows that
$$
     \| Q^B_t f(\Pi_B)
Q^B_s\|= \left\| \left( \psi_s  f \psi_t\right)(\Pi_B) \right\|
     \lesssim\eta(t/s)\sup_{S^o_\mu}|f|$$ for all $ t,s>0,
$ where $\eta(x):=\min\{x,\frac1x\}\left(1+\left|\log|x| \right|
\right)$. A Schur estimate now gives
\begin{align*}
     \|f(\Pi_B)u\|^2 & \approx \int_0^\infty\|Q^B_t
f(\Pi_B)u\|^2\,\frac{dt}t \\
      &\approx\int_0^\infty\left\|\int_0^\infty (Q^B_t
f(\Pi_B)Q^B_s)(Q^B_s
u)\,\frac{ds}s \right\|^2\,\frac{dt}t \\
    & \lesssim\sup_{S^o_\mu}|f|^2 \int_0 ^\infty \left( \int_0^\infty
\eta(t/s)\frac{ds}s \right)
       \left( \int_0^\infty \eta(t/s) \|Q^B_s u\|^2 \frac{ds}s \right)
\frac{dt}{t}\\
    & \lesssim\sup_{S^o_\mu}|f|^2\int_0^\infty\|Q^B_s
u\|^2\,\frac{ds}s\approx\sup_{S^o_\mu}|f|^2\|u\|^2
\end{align*}
for all $u\in\clos{\ran(\Pi_B)}$, which proves that $\Pi_B$ has a
bounded $S^o_\mu$ holomorphic functional calculus in
$L_2(\rnum^n;\cnum^N)$. \qedend \end{proof}

What remains is for us to obtain the estimate (\ref{thesqest})
under all the hypotheses (H1--8). This is achieved in the next
section.

\section{Harmonic analysis of $\Pi_B$}  \label{harmansec}

In this section we prove the square function estimate
(\ref{thesqest}) under the hypotheses (H1--8) stated in Section 2.
By Proposition \ref{whatisneeded}, this then suffices to prove
Theorems \ref{mainthm} and \ref{mainthm2}. This section is an
adaptation of the proof of the Kato square root problem for
divergence-form elliptic operators \cite{HLMc,AHLMcT,AHMcT}, though some estimates require new procedures.  For example, we develop new methods based on hypotheses (H5--6) to prove off-diagonal estimates for resolvents of $\Pi_B$, as the arguments normally used in proving Caccioppoli-type estimates for divergence-form operators do not apply.

We use the following dyadic decomposition of $\rnum^n$. Let
$\dyadic= \bigcup_{j=-\infty}^\infty\dyadic_{2^j}$ where
$\dyadic_t:=\{ 2^j(k+(0,1]^n) :k\in\znum^n \}$ if $2^{j-1}<t\le
2^j$. For a dyadic cube $Q\in\dyadic_{2^j}$, denote by $l(Q)=2^j$
its \emph{sidelength}, and by $R_Q:=Q\times(0,2^j]$ the associated
\emph{Carleson box}. Let the {\em dyadic averaging operator} $A_t:
\mH \longrightarrow \mH $ be given by
$$
       A_t u(x) := u_Q:= \barint_{\hspace{-5pt}Q} u(y)\,  dy = \frac{1}{|Q|} \int_Q
       u(y)\,  dy
$$
for every $x \in \rnum^n$ and $t>0$, where  $Q \in \dyadic_t$ is
the unique dyadic cube  containing $x$.

\begin{defn}
By the {\em principal part} of the operator family $\Theta^B _t$
under consideration, we mean the multiplication operators
$\gamma_t$ defined by
$$
       \gamma_t(x)w:= (\Theta^B _t w)(x)
$$
for every $ w\in \cnum^N$. Here we view $w$ on the right-hand side
of the above equation as the constant function defined on $\rnum
^n$ by $w(x):=w$. It will be proven in Corollary \ref{gammaprops}
that $\gamma_t \in L_2^{\text{loc}}(\rnum^n; \mL(\cnum^N))$.
\end{defn}
To prove the square function estimate (\ref{thesqest}), we
estimate each of the following three terms separately
\begin{equation}
\begin{split} \label{sqfcn2}
       \int_0^\infty&\|\Theta^B _t P_t u\|^2 \frac{dt}t \lesssim
       \int_0^\infty\|\Theta^B _t P_t u-\gamma_t A_t P_t u\|^2\frac{dt}t
\\
       &+ \int_0^\infty\|\gamma_t A_t (P_t-\I) u\|^2\frac{dt}t
     + \int_0^\infty\int_{\rnum^n} |A_t u(x)|^2 |\gamma_t(x)|^2
\frac{dxdt}t
\end{split}
\end{equation}
when $u\in\ran(\Pi)$. We estimate the first two terms in
Section~\ref{princsec}, and the last term in
Section~\ref{Carlesonsec}. In the next section we introduce
crucial off--diagonal estimates for various operators involving
$\Pi_B$, and also prove local $L_2$ estimates for $\gamma_t$.

\subsection{Off--diagonal estimates}

We require off--diagonal estimates for the following class of
operators. Denote $\brac x:=1+|x|$, and
$\dist(E,F) :=\inf\{|x-y|:x\in E,y\in F\}$ for every
$E,F\subset\rnum^n$.

\begin{prop}  \label{pseudoloc}
Let $U_t$ be given by either $R^B_t$ for every nonzero $t \in \R$,
or $P^B_t$, $Q^B_t$ or
       $\Theta^B _t$ for every $t>0$ (see Remark \ref{multiop} and Definition
       \ref{defnope}). Then for every $M \in \nnum$ there exists $C_M>0$
       (that depends only on $M$ and the hypotheses (H1--8))
       such that
\begin{equation} \label{odn}
       \|U_t u\|_{L_2(E)} \le C_M \brac{\dist (E,F)/t}^{-M}\|u\|
\end{equation}
whenever $E,F \subset \R^n$ are Borel sets, and $u \in \mH$
satisfies $\supp u\subset F$.
\end{prop}

\begin{proof}
First consider the resolvents  $R^B_t = (\I+it\Pi_B)^{-1}$
    for all nonzero $t\in\rnum$. As we have already proved uniform bounds
for $R^B_t$ in Proposition \ref{typeomega}, it suffices to prove
$$
       \|(\I+it\Pi_B)^{-1} u\|_{L_2(E)} \le C_M
       (t/\dist(E,F))^M\|u\|
$$
for all disjoint $E$, $F\subset\rnum^n$, for all $|t|\le \dist
(E,F)$, and for all  $u \in \mH$ with $\supp u\subset F$.

We prove this result by induction. Proposition~\ref{typeomega}
proves this statement for $M=0$. Assume that the statement is true
for some given $M\in \nnum$. Write $$\widetilde E:=\{x\in\rnum^n :
\dist(x,E)<\tfrac12\dist(x,F) \}$$ and let
$\eta:\rnum^n\longrightarrow[0,1]$ be a Lipschitz function such
that $\supp\eta\subset\widetilde E$, $\eta|_E=1$ and
$$\|\nabla\eta\|_\infty\leq 4/\dist(E,F).$$ We now use (H5--6) to
calculate that
$$
       [\eta\I,(\I+it\Pi_B)^{-1}]= itR_t^B(\Gamma_{\nabla\eta}+
B_1\Gamma^*_{\nabla\eta}B_2 )R_t^B
$$
and therefore  \begin{equation*} \begin{split} \|(\I+it\Pi_B)^{-1}
u\|_{L_2(E)} & \le \|\eta (\I+it\Pi_B)^{-1}
u\| \\
     &  =\|[\eta\I,(\I+it\Pi_B)^{-1}] u\| \\
      & \lesssim C_0 t\|\nabla\eta\|_{\infty} \|R_t^B u\|_{L_2\left(\widetilde E
      \right)}
      \\ & \lesssim  C_0 t\|\nabla\eta\|_{\infty} C_{M}
(t/\dist(\widetilde E,F))^{M}\|u\| \\ &
       \lesssim C_0 C_{M} (t/\dist(E,F))^{M+1}\|u\|\,.
\end{split} \end{equation*}
This completes the induction step and thus proves the proposition
for the resolvents $R^B_t$.  The result for $P^B_t$ and $Q^B_t$
follows, as they are linear combinations of resolvents.

Now consider $\Theta^B _t= t\Gamma^*_B P^B_t$. We have
$$ \|\Theta^B _t u\|_{L_2(E)} \le  \|\eta\Theta^B _t
u\| \le \|[\eta\I,t\Gamma^*_B]P^B_tu \|+ \|t\Gamma^*_B\eta P^B_t
u\|\,.
$$
By Lemma \ref{lem:gam1} the last term is bounded by
$$\|t\Pi_B\eta P^B_t u\| \le \|[\eta\I,t\Pi_B]P^B_t u \|+\|\eta Q^B_t
u\|$$ and so, using (H6) and the bounds already obtained for
$P^B_t$ and $Q^B_t$, we conclude that for each $M\geq0$,
$$ \|\Theta^B _t
u\|_{L_2(E)} \lesssim   t\|\nabla\eta\|_\infty
\|P^B_tu\|_{L_2\left(\widetilde E \right)} +\|Q^B_t
u\|_{L_2\left(\widetilde E\right)}\lesssim \brac{\dist
(E,F)/t}^{-M}\|u\|\,.$$ This completes the proof. \qedend
\end{proof}

A simple consequence of Proposition \ref{pseudoloc}  is that
\begin{equation}  \label{ODest}
       \|U_s u\|_{L_2(Q)}\le \sum_{R\in\dyadic_t}\|U_s(\chi_Ru)\|_{L_2(Q)}
       \lesssim
       \sum_{R\in\dyadic_t}\brac{\dist(R,Q)/s}^{-M}\|u\|_{L_2(R)}
\end{equation}
whenever $0<s\leq t$ and $Q\in\dyadic_t$, where $U_s$ is as
specified in Proposition \ref{pseudoloc}. We also note that the
dyadic cubes satisfy
\begin{equation}  \label{separationest}
       \sup_{Q\in\dyadic_t}\sum_{R\in\dyadic_t}
\brac{\dist (R,Q)/t}^{-(n+1)}\lesssim 1
\end{equation}
and therefore, choosing $M \ge n+1$, we see that $U_t$ extends to
an operator $U_t:  L_\infty (\rnum^n) \longrightarrow
L_2^{\text{loc}}(\rnum^n)$.

A consequence of the above results with $U_t=\Theta^B_t$ is:

\begin{cor}  \label{gammaprops}
       The functions $\gamma_t\in L_2^{\text{loc}}(\rnum^n;
\mL(\cnum^N))$ satisfy the boundedness conditions
$$
       \barint_{\hspace{-5pt}Q} |\gamma_t(y)|^2 \, d y \lesssim  1$$
for all $Q\in\dyadic_t$.  Moreover $\|\gamma_t A_t\|\lesssim 1$
uniformly for all $t>0$.
\end{cor}

\subsection{Principal part approximation}  \label{princsec}
     In this section we prove the principal part approximation
$\Theta^B _t\approx\gamma_t$ in the sense that we estimate the
first two terms on the right-hand side of (\ref{sqfcn2}). The
following lemma is used in estimating the first term.

\begin{lem} [A weighted Poincar\'e inequality]  \label{poincarelem}
If $Q\in\dyadic_t$ and $\beta < -2n$, then we have
$$
      \int_{\rnum^n}|u(x)-u_Q|^2\brac{\dist(x,Q)/t}^{\beta} \,  dx \lesssim
      \int_{\rnum^n}|t\nabla u(x)|^2 \brac{\dist(x,Q)/t}^{2n+\beta}\,  dx
$$
for every $u$ in the Sobolev space $H^{1} (\rnum^n; \cnum^N)$.
\end{lem}
\begin{proof}
Without loss of generality we may assume that $t=1$ and that $Q$
is the unit cube centred at $x=0$. By \cite[p.\ 164]{truding} we
have
$$ \int_{\R^n} |u(y) - u_Q |^2 \chi_{r}(y) \, dy
\lesssim  \int_{\R^n} |\nabla u(y) | ^2 r^{2n} \chi_{r}(y) \, d y
$$ for every $r \ge 1$, where we write $\chi_r$ to  denote the
characteristic function of $\{ y \in \R^n : |y| \le r\}$.
Integrating the above inequality over $(1,\infty)$ against the
measure $dr^\beta$ gives the desired result. \qedend \end{proof}

We now estimate the first term in the right-hand side of
(\ref{sqfcn2}).

\begin{prop}  \label{secondtermprop}
       For all $u\in \ran(\Pi)$, we have
$$
       \int_0^\infty\|\Theta^B _t P_t u-\gamma_t A_t P_t u\|^2\,\frac{dt}t
       \lesssim  \|u\|^2.
$$
\end{prop}
\begin{proof}
Using Proposition~\ref{pseudoloc}, estimate (\ref{separationest}),
Lemma~\ref{poincarelem} and then the coercivity assumption (H8),
we get for any $v\in \ran(\Pi)$, that
\begin{equation*} \begin{split}
       \|\Theta^B _t v-\gamma_t A_t v  \|^2
      & =\sum_{Q\in\dyadic_t}\| \Theta^B _t(v-v_Q) \|_{L_2(Q)}^2 \\
       & \lesssim \sum_{Q\in\dyadic_t}\Big(
\sum_{R\in\dyadic_t}\brac{d(R,Q)/t}^{-(3n+1)}\|v-v_Q\|_{L_2(R)}\Big)^2
      \\ & \lesssim \sum_{Q\in\dyadic_t}
       \int_{\rnum^n}|v(x)-v_Q|^2\brac{d(x,Q)/t}^{-(3n+1)} \\
      & \lesssim \sum_{Q\in\dyadic_t}\int_{\rnum^n}|t\nabla v(x)|^2
\brac{d(x,Q)/t}^{-(n+1)}
      \lesssim \|t\nabla v \|^2 \lesssim \|t\Pi v\|^2
\end{split} \end{equation*}
and therefore, taking $v=P_tu$ and using (\ref{standardsqest}),
that
$$
       \int_0^\infty\|\Theta^B _t P_t u-\gamma_t A_t P_t u\|^2\,\frac{dt}t
       \lesssim \int_0^\infty\|Q_t u\|^2\,\frac{dt}t
       \lesssim  \|u\|^2.\quad  \qedend
$$
 \end{proof}

We use the following lemma to estimate  the second term in the
right-hand side of (\ref{sqfcn2}), and also in the proofs of
Lemmas \ref{lem:ncar1} and \ref{lem:ncar10}. (c.f. Lemma 5.15 of
\cite{AHLMcT}.)

\begin{lem} \label{interpolationlemma}
Let $\Upsilon$ be either $\Pi$, $\Gamma$ or $\Gamma^*$. Then we
have the estimate
\begin{equation}  \label{interpolest}
\left| \barint_{\hspace{-5pt} Q} \Upsilon u  \right|^2 \lesssim
\frac 1{l(Q)} \left( \barint_{\hspace{-5pt}Q} |u|^2 \right)^{1/2}
\left( \barint_{\hspace{-5pt}Q} | \Upsilon u|^2 \right)^{1/2}
\end{equation}
for all $Q \in\dyadic$ and $u \in \dom(\Upsilon)$.
\end{lem}
\begin{proof}
Let $t=(\int_Q|u|^2)^{1/2}(\int_Q|\Upsilon u|^2)^{-1/2}$. If $t\ge
\frac14 l(Q)$, then (\ref{interpolest}) follows directly from the
Cauchy--Schwarz inequality. If $t\leq \frac14 l(Q)$, let $\eta\in
C^\infty_0(Q)$ be a real-valued bump function such that
$\eta(x)=1$ when $\dist(x,\rnum^n\setminus Q)>t$, and
$|\nabla\eta|\lesssim 1 /t$. Using the cancellation property (H7)
of $\Upsilon$ and the Cauchy--Schwarz inequality, we obtain
\begin{equation*} \begin{split}
       \left|\int_Q \Upsilon u\right| &=\left|\int_Q \eta \Upsilon u+
\int_Q (1-\eta)\Upsilon u\right|
        =\left|\int_Q [\eta,\Upsilon ]u+ \int_Q (1-\eta)\Upsilon u\right| \\
       &\lesssim  \|\nabla\eta\|_{\infty} ( t l(Q)^{n-1})^{1/2}
\left(\int_Q|u|^2\right)^{1/2} + (t l(Q)^{n-1})^{1/2}
\left(\int_Q|\Upsilon u|^2\right)^{1/2}
\end{split} \end{equation*}
which gives (\ref{interpolest}) on substituting the chosen value
of $t$. \qedend \end{proof}

We now estimate the second  term in the right-hand side of
(\ref{sqfcn2}).

\begin{prop}  \label{thirdtermprop}
       For all $u\in \mH$, we have
$$
       \int_0^\infty \|\gamma_t A_t (P_t-\I) u\|^2\,\frac{dt}t \lesssim  \|u\|^2.
$$
\end{prop}
\begin{proof}
Corollary~\ref{gammaprops} shows that $\|\gamma_t A_t\|\lesssim 1$
and since $A_t^2=A_t$ it suffices to prove the square function
estimate with integrand $\|A_t (P_t-\I) u\|^2$. If $u\in\nul(\Pi)$
then this is zero. If $u\in\clos{\ran(\Pi)}$ then write
$u=2\int_0^\infty Q_s^2u\frac{ds}s$. The result will follow from
another Schur estimate and (\ref{standardsqest}) once we have
obtained the bound
$$
       \|A_t (P_t-\I)Q_s\|\lesssim \min\{\tfrac st,\tfrac ts\}^{1/2}$$
for all $s,t >0$.

Note that $  (\I-P_t)Q_s = \tfrac ts Q_t(\I-P_s)$
 and $ P_tQ_s =\tfrac st Q_tP_s$
for every $s,t>0$. Thus, if $t \le s$, then $$\|A_t
(P_t-\I)Q_s\|\lesssim \|(P_t-\I)Q_s\|\lesssim  t/s,$$ while if
$t>s$, then $$\|A_t (P_t-\I)Q_s\|\lesssim
\|P_tQ_s\|+\|A_tQ_s\|\lesssim s/t+\|A_tQ_s\|.$$ To estimate
$\|A_tQ_s\|$, we use Lemma~\ref{interpolationlemma} with
(\ref{ODest}) and (\ref{separationest}) to obtain
\begin{equation*} \begin{split}
       \|A_tQ_s u\|^2 &= \sum_{Q\in\dyadic_t} |Q|\,
       \bigg| \barint_{\hspace{-5pt}Q} s\Pi(\I+s^2 \Pi ^2 )^{-1}u \bigg|^2 \\
       & \lesssim \frac st\sum_{Q\in\dyadic_t}
       \bigg( \int_{Q}|P_su|^2 \bigg)^{1/2}
       \bigg( \int_{Q}|Q_su|^2 \bigg)^{1/2} \\
       & \lesssim \frac st\sum_{Q\in\dyadic_t}
         \bigg(
\sum_{R\in\dyadic_t}\brac{d(R,Q)/t}^{-(n+1)}\|u\|_{L_2(R)}
\bigg)^2
     \\
       & \lesssim \frac st \sum_{Q\in\dyadic_t}
         \bigg( \sum_{R'\in\dyadic_t}\brac{d(R',Q)/t}^{-(n+1)} \bigg)
         \bigg( \sum_{R\in\dyadic_t}\brac{d(R,Q)/t}^{-(n+1)}\|u\|_{L_2(R)}^2
         \bigg)\\& \lesssim \frac st \|u\|^2
\end{split} \end{equation*}
which completes the proof. \qedend \end{proof}

We have now estimated the first two terms in the right-hand side
of (\ref{sqfcn2}).

\subsection{Carleson measure estimate}  \label{Carlesonsec}

In this subsection we estimate the third term in the right-hand side
of (\ref{sqfcn2}). To do this we reduce the problem to a Carleson
measure estimate, drawing upon the ``$T(b)$" procedure developed by Auscher and Tchamitchian  \cite[Chapter 3]{ATch}. Recall that a measure $\mu$ on $\R^n \times
\R^+$ is said to be \emph{Carleson}  if
$\|\mu\|_{\mC}:=\sup_{Q\in\triangle}|Q|^{-1}\mu(R_Q)<\infty$.
 Here
and below $R_Q := Q \times (0,l(Q)]$ denotes the {\it Carleson
box} of any cube $Q$.
 We recall the
following theorem of Carleson.

\begin{thm}\cite[p.\ 59]{stein:harm}
If $\mu$ is a Carleson measure on $\R^n \times \R^+$ then
$$ \iint_{\R^n \times (0,\infty)} |A_t u (x)|^2 \, d \mu(x,t) \le C
\|\mu \|_{\mC} \| u \|^2$$ for every $u \in \mH$. Here $C>0$ is a
constant that depends only on $n$.
\end{thm}

Thus, in order to prove (\ref{sqfcn2}) it suffices to show that
\begin{equation} \label{ncar1} \iint_{R_Q} | \gamma_t (x) |^2 \,
\frac{dxdt}{t} \lesssim |Q|\end{equation} for every dyadic cube $Q
\in \triangle$. Following \cite{AHLMcT} or more precisely
\cite{AHMcT}, we set $\sigma>0$; the exact value to be chosen
later. Let $\mV$ be a finite set consisting of $\nu \in  \mL
(\C^N)$ with $|\nu|=1$, such that $\bigcup _{ \nu \in \mV} K_\nu =
\mL ( \C^N) \setminus \{0\},$ where
$$K_\nu := \left\{ \nu' \in \mL(\C^N) \setminus \{0\} : \left| \frac{\nu'}{|\nu'|} - \nu \right|
\le \sigma  \right\}.$$ To prove (\ref{ncar1}) it suffices to show
that
\begin{equation} \label{ncar2}
\iint_{\substack{ (x,t) \in R_Q \\ \gamma_t(x) \in K_\nu }} |
\gamma_t(x) |^2 \, \frac{dxdt}{t} \lesssim |Q|\end{equation} for
every $\nu \in \mV$. By the John-Nirenberg lemma for Carleson measures as applied
in \cite[Section 5]{AHLMcT}, in order to prove (\ref{ncar2}) it
suffices to prove the following claim.

\begin{prop} \label{prop:ncar2}
There exists $\beta
>0$ such that for every dyadic cube $Q\in \triangle$ and $\nu \in \mL
( \C^N)$ with $|\nu |=1$, there is a collection $\{Q_k\}_{k}
\subset \triangle$ of disjoint subcubes of $Q$ such that $
|E_{Q,\nu}|
>  \beta |Q|$ where $E_{Q,\nu }= Q \setminus \bigcup_k Q_k$, and such
that
$$ \iint_{\substack { (x,t) \in  E^* _{Q,\nu} \\ \gamma_t(x) \in
K_\nu }} | \gamma_t(x) |^2 \,\frac{dx dt}{t} \lesssim |Q|$$ where
$E^* _{Q ,\nu } = R_Q \setminus \bigcup_k R_{Q_k}$.
\end{prop}

Fix a dyadic cube $Q \in \triangle$ and fix $\nu \in \mL(\C^N)$
with $| \nu |=1$. Choose $\hat w, w \in \C^N$ with $|\hat
w|=|w|=1$ and $\nu^* (\hat w )=w$. Let $\eta_Q$ be a smooth cutoff
function with range $[0,1]$, equal to $1$ on $2Q$, with support in
$4Q$, and such that $\| \nabla \eta_Q \|_\infty \le \tfrac{1}{l}$
where $l=l(Q)$. Define $w_Q := \eta_Q w $, and for each
$\epsilon>0$, let
$$ f_{_{Q,\epsilon}}^w :=w_{_Q}-\epsilon l
i\Gamma(1+\epsilon l i\Pi_B)^{-1}w_{_Q}\\
=\left(1+\epsilon l i\Gamma_B^*\right)(1+\epsilon l
i\Pi_B)^{-1}w_{_Q}\,.$$

\begin{lem} \label{lem:ncar1}
We have $\|f^w _{Q,\epsilon } \| \lesssim |Q|^{1/2},$
\begin{equation*}
\iint_{R_Q} |\Theta_t ^B f^w_{Q,\epsilon}|^2 \, \frac{dxdt}{t}
\lesssim \frac{1}{\epsilon^2} |Q| \quad \text{and}\end{equation*}
\begin{equation*} \left| \barint_{\hspace{-5pt}Q} f^w_{Q,\epsilon} - w \right| \le c
\epsilon^{1/2}\end{equation*} for every $\epsilon >0$. Here $c>0$
is a constant that depends only on  hypotheses (H1--8).
\end{lem}

\begin{proof}
The first estimate  can be deduced from Proposition
\ref{typeomega} and Lemma \ref{lem:gam1}. To obtain the second
estimate, observe by the nilpotency of $\Gamma^*_B$ that
\begin{equation*}
\begin{split} \Theta_t ^B f^w_{Q,\epsilon} &= (\I +
t^2 \Pi_B^2)^{-1} t \Gamma^*_B (\I+ \epsilon l i \Gamma^*_B) (\I +
\epsilon l i \Pi_B)^{-1}w_Q \\
&= \tfrac{t}{\epsilon l} (\I + t^2 \Pi_B^2)^{-1} \epsilon l
\Gamma^*_B  (\I + \epsilon l i \Pi_B)^{-1}w_Q
\end{split} \end{equation*}
and therefore by Proposition \ref{typeomega} and Lemma
\ref{lem:gam1} that
$$\iint_{R_Q}
|\Theta_t ^B f^w_{Q,\epsilon}|^2 \, \frac{dxdt}{t} \lesssim |Q|
\int_{0} ^l
   \left(\frac{t}{\epsilon l}\right)^2 \, \frac{dt}{t} \lesssim
  \frac{1}{\epsilon^2} |Q|.
$$
To obtain the last estimate, we use Lemma~\ref{interpolationlemma}
with $\Upsilon=\Gamma$ and $u=(\I+\epsilon li\Pi_B)^{-1}w_Q$ to
show that
\begin{multline*}
       \left| \barint_{\hspace{-5pt}Q} f^w_{Q,\epsilon} - w\right|
       = \left|\barint_{\hspace{-5pt}Q} \epsilon l\Gamma (\I+\epsilon li\Pi_B)^{-1}w_Q\right|
     \\
       \lesssim \epsilon^{1/2} \left(\barint_{\hspace{-5pt}Q} |(\I+\epsilon
li\Pi_B)^{-1}w_Q|^2\right)^{1/4}
       \left(\barint_{\hspace{-5pt}Q} | \epsilon l \Gamma (\I+\epsilon
li\Pi_B)^{-1}w_Q|^2\right)^{1/4}
       \lesssim  \epsilon^{1/2}.
\end{multline*}
This completes the proof. \qedend \end{proof}

For the choice $\epsilon = \tfrac{1}{4c^2}$, let $f^w_Q
=f^w_{Q,\epsilon}.$ The above lemma implies that
$$ \re \left(w, \barint_{\hspace{-5pt}Q} f^w_Q\right ) \ge \frac{1}{2}\,.
$$

\begin{lem}
There exists $\beta,c_1,c_2 >0$ that depend only on (H1--8), and
there exists a collection $\{Q_k\}$ of dyadic subcubes of $Q$ such
that $ |E_{Q,\nu}| > \beta |Q|$ where $E_{Q,\nu }= Q \setminus
\bigcup_k Q_k$, and such that \begin{equation} \label{recar12} \re
\left( w , \barint_{\hspace{-5pt}Q'} f^w_{Q} \right) \ge c_1 \quad
\text{and} \quad \barint_{\hspace{-5pt}Q'} | f^w_Q| \le c_2
\end{equation} for all dyadic subcubes $Q' \in \triangle$ of $Q$
which satisfy $R_{Q'} \cap E^* _{Q,\nu} \neq \emptyset,$ where
$E^* _{Q, \nu } = R_Q \setminus \bigcup_k R_{Q_k}$.
\end{lem}

\begin{proof}
Fix $\alpha >0$. Let $\CB_1 \subset \triangle$ be the collection
of maximal dyadic subcubes $S \in \triangle$ of $ Q$ such that
$$ \re\left (w, \barint_{\hspace{-5pt}S} f^w _Q\right ) < \alpha$$ and let
$\CB_2 \subset \triangle$ be the collection of maximal  dyadic
subcubes  $S \in \triangle $ of $ Q$ such that
$$\barint_{\hspace{-5pt}S} | f^w _Q|
>  \frac{1}{\alpha}.$$
Let $\{Q_k\}$ be an enumeration of the maximal cubes in $\CB_1
\cup \CB_2$. These are the bad cubes. By construction we have each
dyadic subcube $Q' \in \triangle$ of $Q$ with $R_{Q'} \cap E^*
_{Q,\nu} \neq \emptyset$  satisfies (\ref{recar12}) with
$c_1=\alpha$ and $c_2=\tfrac{1}{\alpha}$. These are the good
cubes. Thus, to prove the lemma it suffices to show that for an
appropriate choice of $\alpha>0$, that depends only on (H1--8),
there exists $\beta>0$ such that $|E_{Q,\nu}| > \beta |Q|$.

We  use the rough estimate
$$ |E_{Q,\nu}| \ge | Q \setminus \bigcup \CB_1| - | \bigcup \CB_2|.$$
By construction and by Lemma \ref{lem:ncar1} we have
$$ | \bigcup \CB_2 | = \sum_{S \in \CB_2} |S| \le \alpha
\sum_{S \in \CB_2} \int_S |f^w _Q| \le \alpha \int_Q |f^w _Q|
\lesssim \alpha |Q|$$ and
\begin{equation*}
\begin{split} \frac{1}{2}|Q| \le \re\left (w, \int_Q f^w _Q\right) &=
\sum_{S \in \CB_1} \re\left (w, \int_S f^w _Q\right ) + \re\left
(w, \int_{Q \setminus \bigcup \CB_1} f^w _Q\right )
\\ &\lesssim \alpha \sum_{S \in \CB_1}  |S| + \left( \int_Q |f^w _Q |^2
\right)^{1/2} |Q \setminus \bigcup \CB_1|^{1/2} \\ &\lesssim
\alpha |Q| + |Q|^{1/2}
  | Q \setminus \bigcup \CB_1 |^{1/2}. \end{split} \end{equation*}
The desired estimate follows by a sufficiently small choice of
$\alpha>0$ that depends  only on (H1--8). This completes the
proof. \qedend \end{proof}

We now choose $\sigma = \tfrac{c_1}{2c_2}$.

\begin{lem} \label{lem:ncar10}
If $(x,t) \in E^*_{Q,\nu}$ and $\gamma_t(x) \in K_\nu$ then
$$ \left| \gamma_t(x) \left(A_t f^w_Q (x) \right) \right|
  \ge \tfrac{1}{2} c_1|\gamma_t(x)| .$$
\end{lem}

\begin{proof}
To see the result apply the previous lemma to deduce that $$\left|
\nu \left( A_t f^w _Q (x) \right) \right| \ge \re \left(\hat w ,
\nu \left( A_t f^w_Q (x ) \right) \right) = \re \left( w, A_t f^w
_Q (x) \right) \ge c_1$$ and then furthermore that
$$\left| \frac{\gamma_t(x)}{|\gamma_t(x)|} \left(A_t f^w _Q (x)
\right) \right| \ge \left| \nu \left(A_tf^w _Q (x) \right) \right|
- \left| \frac{\gamma_t(x)}{|\gamma_t(x)|} - \nu \right| \left|
A_t f^w _Q (x) \right| \ge c_1 - \sigma  c_2 = \tfrac{1}{2} c_1.$$
\qedend \end{proof}

\begin{proof} [Proof of Proposition \ref{prop:ncar2}]

By Lemma  \ref{lem:ncar10} we have
\begin{equation*} \begin{split} \iint_{\substack { (x,t) \in  E^*
_{Q,\nu} \\ \gamma_t(x) \in K_\nu }} | \gamma_t(x) |^2 \,\frac{dx
dt}{t} & \lesssim \iint_{R_Q} \left|\gamma_t(x) \left(A_t f^w _Q
(x) \right) \right|^2 \,\frac{dx dt}{t} \\& \lesssim \iint_{R_Q}
\left| \Theta_t ^B f^w_Q - \gamma_t A_t f^w_Q  \right|^2 \,
\frac{dxdt}{t} + \iint_{R_Q} \left|\Theta_t ^B f^w_Q \right|^2 \,
\frac{dxdt}{t}\,.\end{split}
\end{equation*}
Lemma \ref{lem:ncar1} implies that the last term in the above
inequality is bounded by a constant (that depends only on (H1--8))
times $|Q|$.

It remains to show that
\begin{equation} \label{eq:car21}\iint_{R_Q} \left| \Theta_t ^B
f^w_Q - \gamma_t A_t f^w_Q  \right|^2 \, \frac{dxdt}{t} \lesssim
|Q|.\end{equation}
Observe that
\begin{equation} \label{eq:car32} \Theta_t^B f^w_Q - \gamma_t A_t f^w_Q
  = - \left( \Theta_t^B - \gamma_t  A_t \right)  \epsilon l i
\Gamma (1 + \epsilon l i \Pi_B)^{-1} w_Q + (\Theta _t ^B
  - \gamma_t  A_t) w_Q \,.
\end{equation}
Since $ \epsilon l i \Gamma (1 + \epsilon l i \Pi_B)^{-1} w_Q \in
\ran (\Gamma)$, we have by the results of Sections \ref{opthsec}
and \ref{princsec} (specifically, part (i) in the proof of
Proposition \ref{whatisneeded}, and also \ref{secondtermprop} and
\ref{thirdtermprop}) that
$$ \iint_{R_Q} \left| \left( \Theta_t^B - \gamma_t  A_t \right)  \epsilon l i
\Gamma (1 + \epsilon l i \Pi_B)^{-1} w_Q \right|^2 \,
\frac{dxdt}{t} \lesssim \|w_Q \|^2 \lesssim |Q|.$$ We also have
$$
       (\Theta^B _t-\gamma_t A_t) w_Q(x)= \Theta^B _t((\eta_Q-1)w)(x)
$$
for every $x \in Q$ and $t>0$. Since $\left(\supp(\eta_Q-1)w
\right) \cap 2Q=\emptyset,$ then (\ref{ODest}) implies that
$$\int_Q|\Theta^B _t((\eta_Q-1)w)(x)|^2 \, dx \lesssim  \frac{t|Q|}{l}$$
when $0<t\leq l$, and therefore that
$$
\iint_{R_Q} \left|   (\Theta^B _t-\gamma_t A_t) w_Q(x) \right|^2
\, \frac{dxdt}{t} \lesssim |Q|.
$$
This proves (\ref{eq:car21}) and so completes the proof of
Proposition \ref{prop:ncar2}. \qedend \end{proof}

\begin{proof}[Proof of Theorems \ref{mainthm} and \ref{mainthm2}]
We have demonstrated in this section that the square function
estimate (\ref{thesqest}) holds for all $u\in\ran(\Pi)$ and some
constant $c$ which depends only on the bounds in (H1--8). These
hypotheses are invariant on replacing $\{\Gamma,B_1, B_2\}$ by
$\{\Gamma^* ,B_2, B_1\}$, $\{\Gamma^*,B_2^*, B_1^*\}$ and
$\{\Gamma,B_1^*, B_2^*\}$. So, by Proposition \ref{whatisneeded},
we conclude that $\Pi_B$ satisfies the quadratic estimate
(\ref{eq:quad}) for all $u\in\overline{\ran(\Pi_B)}$, and has a
bounded holomorphic $S^o_\mu$ functional calculus. \qedend
\end{proof}

\section{Holomorphic dependence} \label{sect:holo}

In this section we show that under the appropriate hypotheses,
resolvents, projections, bounded members of the functional
calculus, and quadratic estimates, all depend holomorphically on
holomorphic perturbations of $\Pi_B$. Recall that if $\mH$ and
$\mK$ are Hilbert spaces and $U \subset \C$ is open, then an
operator valued function $T: U \to \mL(\mH,\mK)$ is said to be
\emph{holomorphic} if it is (complex) differentiable in the
uniform topology everywhere in $U$.

\begin{thm} \label{new1}
Let $U \subset \cnum$ be open,  let $B_1, B_2: U \to \mL(\mH)$ be
holomorphic functions such that $B_1(z)$ and $ B_2(z)$ satisfy
(H1--3) uniformly for each $z\in U$, and let $\tau \in \cnum
\setminus S_\mu$. Then
    the function given by $z \mapsto
(1+\tau \Pi _{B(z)})^{-1}$ is holomorphic on $U$, the function
given by $z \mapsto {\mathbf P}^0 _{B(z)}$ is holomorphic on $U$,
and the function given by $z \mapsto \psi(\Pi_{B(z)})$ is
holomorphic on $U$ for every $\psi \in \Psi(S^o _ \mu)$.
\end{thm}

\begin{rem} \label{remobs}
An interesting observation that arose from our consideration of
Theorem \ref{new1} is that under its hypotheses, not only is the
function given by $z \mapsto \mathbf{P}^0 _{B(z)}$ holomorphic on
$U$, but so too are the functions given by $z \mapsto \mathbf{P}^1
_{B(z)}$ and $z \mapsto \mathbf{P}^2 _{B(z)}$. This means that the
Hodge decomposition (\ref{Hodge}) is holomorphic on $U$.

Moreover, we have
\begin{equation} \label{eq:deri1} \frac{d}{dz} \PP^0  =  -\PP ^0 A_1
\Pa - \Pb  A_2 \PP^0, \end{equation}
\begin{equation} \label{eq:deri2} \frac{d}{dz} \PP^1 =  \PP ^0 A_1
\Pa - \Pb A_2 \PP^1, \end{equation} and
\begin{equation} \label{eq:deri3} \frac{d}{dz} \PP^2 =  -\PP ^2 A_1
\Pa + \Pb A_2 \PP^0 \, .\end{equation} Here $A_1(z) =
\frac{d}{dz}B_1(z)$ and $A_2(z) = \frac{d}{dz}B_2(z)$, and the
operators $\Pa$ and $\Pb$ in $\mL(\mH)$ are defined in Appendix
\ref{append}, and satisfy $\PP^1=B_1\Pa$ and $\PP^2=\Pb B_2$. The
claims of this remark are verified in the Appendix \ref{append}.
\end{rem}

Before proving Theorem \ref{new1} we recall some standard results
from operator theory.  The function $T: U \to \mL(\mH,\mK)$ is
holomorphic if and only if it is locally uniformly bounded (that
is, uniformly bounded on each  compact subset of $U$), and
strongly differentiable (see \cite[p.\ 365]{Kato}). Cauchy's
Theorem, and indeed many standard results about complex-valued
holomorphic functions extend to the operator valued setting. A
suitable reference is \cite[III.14]{DS}. In particular, the
following holds:

\begin{lem} \label{lem:oper}
Let $U \subset \cnum$ be an open set, and let $T_n,T : U
\longrightarrow \mL(\mH,\mK)$ be functions with $T_n$ holomorphic
      for each $n \in \nnum$. Suppose that $T_n(z) u \to T(z) u$ as $n
\to \infty$,
      for every
      $z \in U$ and $u \in \mH$, and that  for every
compact $K \subset U$ there exists  $L>0$ such that $\|T_n (z)\|
\le L$ for every $z \in K$ and $n \in \nnum$. Then $T$ is
holomorphic, and moreover for every $u \in \mH$,  we have  $(T_n
u)$ and $(\frac{d}{dz} T_n u)$ converge locally uniformly to $Tu$
and $\frac{d}{dz} T u$ respectively. (i.e.  the convergence is
uniform on each compact subset of $U$.)
\end{lem}

A sequence $(T_n) \subset \mL(\mH)$  is said to converge  to $T
\in \mL(\mH)$ \emph{strongly}   if for every $u \in \mH$ we have
$\| T_n u - T u \| \to 0$ as $n \to \infty$. We use the fact that,
for any pair of sequences $(S_n), (T_n) \subset \mL(\mH)$ with
$S_n \to S$ and $T_n \to T$ strongly as $n \to \infty$, where $S,T
\in \mL(\mH)$, then $S_n T_n \to S T$ strongly.

\begin{proof}[Proof of Theorem \ref{new1}]
Fix $\tau \in \C \setminus S^o _\mu$. Then
\begin{equation} \label{neq1}
\begin{split} \frac{d}{dz} (\I + \tau \Pi_B)^{-1} &= -(\I + \tau
\Pi_B)^{-1} A_1 \tau \Gamma^* B_2 (\I + \tau \Pi_B)^{-1} \\
&\quad- (\I + \tau \Pi_B)^{-1}  B_1\tau \Gamma^* A_2 (\I + \tau
\Pi_B)^{-1}\end{split}
\end{equation}
where $A_1(z) = \frac{d}{dz}B_1(z)$ and $A_2(z) =
\frac{d}{dz}B_2(z)$. The fact that the above operators are all
uniformly bounded can be obtained from (\ref{equiv1}),
(\ref{equiv2}) and Lemma \ref{lem:gam1}. This proves the first
claim. Thus
    $\{ z \mapsto (\I + in \Pi_B)^{-1}\}_n$ is a collection of
    uniformly bounded functions holomorphic on $U$. Moreover $ \PP^0u =
\lim_{n \to \infty} (\I + in \Pi_{B(z)})^{-1}u$ for all $u\in
\mH$. (This is proved in a setting similar to ours in
\cite[Theorem 3.8]{cowl}; we also prove it as a part of Lemma
\ref{lem:multi}). The second claim now follows from Lemma
\ref{lem:oper}. We now prove the third claim. Fix $\psi \in
\Psi(S_\mu ^o)$. The desired result can now be deduced from the
first claim and Lemma \ref{lem:oper} by using a Riemann sum to
approximate the contour integral representation of
$\psi(\Pi_{B(z)})$ as in (\ref{dunfordschwarz}). This completes
the proof. \qedend \end{proof}

We now adopt the notation from hypotheses (H1--8) and consider the
Hilbert space
$$\mK = L_2 \left( \rnum^n \times (0,\infty), \frac{dx dt}{t} ;
\cnum^N  \right)$$ and for every $\psi \in \Psi (S_\mu ^o)$ and $z
\in U$, define the operator $S_{B(z)}(\psi) : \mH \longrightarrow
\mK$ by
$$ (S_{B(z)}(\psi) u)(x,t) = \left( \psi (t \Pi_{B(z)}) u \right)(x)$$
for every $ u \in \mH$,  $t >0$ and almost every $x \in \rnum^n$.

\begin{thm} \label{newthm2}
Let $U \subset \cnum$ be open,  let $B_1, B_2: U \to \mL(\mH)$ be
holomorphic functions such that $B_1(z)$ and $ B_2(z)$ satisfy
(H1--8) uniformly for each $z\in U$, and let $\omega < \mu <
\frac\pi2$. Then the function given by $z \mapsto f(\Pi_{B(z)})$
is holomorphic on $U$ for every bounded $f: S^o _\mu \cup \{ 0 \}
\longrightarrow \C$ holomorphic on $S^o _\mu$,  and the function
given by $z \mapsto S_{B(z)}(\psi)$ is holomorphic on $U$ for
every $\psi \in \Psi(S^o _ \mu)$.
\end{thm}

\begin{proof}
We prove the first claim. Let $f $ be as above. Since by Theorem
\ref{new1}, the function  $z \mapsto \mathbf{P}^0 _{B(z)}$ is
holomorphic on $U$, we can without loss of generality further
assume that $f(0)=0$. Choose a uniformly bounded sequence
$(\psi_n) \subset \Psi(S_\mu ^o)$ that converges locally uniformly
to  $f$ on $S_\mu ^o$. By Theorem \ref{new1} we have each function
$z \mapsto \psi_n(\Pi_{B(z)})$ is holomorphic on $U$. Moreover, by
Theorem \ref{mainthm2} and (\ref{evenmorefun}), we have that
$\psi_n(\Pi_{B(z)})$ is uniformly bounded (with respect to $n \in
\N$ and $z \in U$) and that $\left(\psi_n(\Pi_{B(z)})\right)$
  converges strongly to $f(\Pi_{B(z)})$ for every $z \in U$. The
first claim of Theorem \ref{newthm2} now follows from Lemma
\ref{lem:oper}.

We now prove the second claim. Let $n \in \nnum$, and define
$\psi^n _t: S_\mu ^o \longrightarrow \cnum$ by  $\psi^n _t (\zeta)
= \psi (t\zeta)$ whenever $\zeta \in S_\mu ^o$ and $1/n < t < n$,
and $\psi^n _t = 0$ otherwise. Next let $S_{B(z)} ^n(\psi)  : \mH
\longrightarrow \mK$ be given by
$$ (S_{B(z)} ^n(\psi)  u)(x,t) = \left( \psi_t ^n (\Pi_{B(z)}) u \right)(x)$$
for every $z \in U$, $ u \in \mH$,   $t >0$ and almost every $x
\in \rnum^n$. We deduce from Theorem \ref{new1} that for every
$t>0$, the function $z \mapsto \psi_t ^n ( \Pi_{B(z)})$ is
holomorphic on $U$, and by Theorem \ref{mainthm2} that this family
of  functions  is uniformly bounded with respect to $t>0$. This
with the fact that $\psi^n _t$ is only non-zero for $t \in
(1/n,n)$ allows us to deduce that the function given by $z \mapsto
S_{B(z)} ^n(\psi) $ is holomorphic on $U$. However, by Remark
\ref{otherpsi} we have $\| S_{B(z)} ^n(\psi) \|$ is uniformly
bounded over every $z \in U$ and $n \in \nnum$, and that $S_{B(z)}
^n(\psi)$ strongly converges to $S_{B(z)}(\psi) $ as $n \to
\infty$ for every $z \in U$. The second claim now follows from
Lemma \ref{lem:oper}. This completes the proof. \qedend
\end{proof}

We use the previous theorem to prove Lipschitz estimates on
members of the functional calculus of the perturbed Dirac operator
$\Pi_B$, and Lipschitz estimates on quadratic functions of
$\Pi_B$.

\begin{thm} \label{prop:lips}
Let $\mH, \Gamma, B_1,B_2,\kappa_1,\kappa_2$ and $n$ be as
outlined in (H1--8). For $i=1,2$, fix $\eta_i < \kappa_i,$ and
then let $ 0< \hat \omega_i <\frac\pi2$ be given by $\cos \hat
\omega_i = \frac{\kappa_i - \eta_i}{\|B_i  \| + \eta_i } .$ Next
let $\hat \omega = \frac12( \hat \omega_1 + \hat \omega_2)$ and
$\hat \omega  < \mu < \frac\pi2$. Then we have
$$ \| f(\Pi_{B}) -f (\Pi
_{B+A}) \| \lesssim ( \|A_1 \|_\infty + \|A_2 \|_\infty) \|f
\|_\infty   $$ for every bounded $f : S_\mu ^o \cup \{0 \}
\longrightarrow \cnum$ holomorphic on $S_\mu ^o$, and every
$A_i\in L_\infty(\rnum^n,\mL(\cnum^N))$ with $\| A_i \|_\infty \le
\eta_i$. Moreover, given $\psi \in \Psi(S_\mu ^o)$, we have
\begin{equation*} \int_0 ^\infty \| \psi(t \Pi_{B}) u- \psi(t
\Pi_{B+A})u \|^2 \frac{dt}{t} \lesssim ( \|A_1\|_\infty^2 + \| A_2
\|_\infty^2) \|u\|^2 \end{equation*}  for all $u \in \mH$, and
every $A_i\in L_\infty(\rnum^n,\mL(\cnum^N))$ with $\| A_i
\|_\infty \le \eta_i$.
\end{thm}

\begin{proof}
For each $i=1,2$, define the functions $G_i:\cnum \longrightarrow
\mL (\mH)$ by $z \mapsto B_i + z A_i,$ and let $$U = \left\{ z \in
\cnum : |z| \le \min \left\{ \eta_1 \|A_1 \|^{-1} , \eta_2 \|A_2
\|^{-1} \right\} \right\}.$$ For all $z \in U$ and  $i=1,2$ we
have
$$ \re ((B_i + z A_i) u,u) \ge (\kappa_i -\eta_i) \|u\|^2$$
for every $u \in \mH$, and therefore
$$  \cos \sup_{u\in\ran(\Gamma^*)\setminus\{0\}}|\arg((B_i +z A_i)
u,u)| \ge   \frac{\kappa_i - \eta_i}{\|B_i  \| + \eta_i } = \cos
\hat \omega_i .$$ We conclude that
     $G_1(z)$ and $G_2(z)$ satisfy (H2) with $\omega_1$ and
$\omega_2$ replaced by $\hat \omega_1$ and $\hat \omega_2$, and
thence by Theorem \ref{newthm2}, that the function given by $z
\mapsto \Pi_{G(z)}$ is holomorphic on $U$. The first claim of the
theorem then follows by Schwarz's Lemma. The second claim is
proved by a similar argument. \qedend \end{proof}

\section{Applications to Riemannian manifolds} \label{sect:riem}

We now consider applications to compact Riemannian manifolds $M$
with metric $g$. For each $x \in M$  let $\wedge T^*_x M$ denote
the complex exterior algebra over the cotangent space $T^*_x M$.
We then let $\wedge T^* M$ and $\mL_M$ denote the bundles over $M$
whose fibres at each $x \in M$ are given by $\wedge T^*_x M$ and
$\mL ( \wedge T^* _x M)$, respectively. We let $\mH= L_2(\wedge
T^* M)$ denote the collection of $L_2$ integrable sections of
$\wedge T^* M$, and let $L_\infty(\mL_M)$  denote the bounded
measurable sections of  $\mL_M$. We let $d^{*} _g$ denote the dual
of $d$ in $\mH$, and consider the Hodge--Dirac operator $D_g:= d +
d^*_g$.

\begin{thm} \label{thm:man1}
Let $M$ be a compact Riemannian manifold with metric $g$, let $B
\in L_\infty(\mL_M)$ be invertible and so that there exists
$\kappa>0$ such that for almost every $x \in \rnum^n$, we have
$$ \re(B(x) v,v) \ge
\kappa | v |^2$$ for every $v \in \wedge T^* _x M$. Let  $ \omega
< \mu < \frac\pi2$ where
$$ \omega := \esssup_{\substack{ x \in M \\ v \in \wedge T_x ^* M}}
  |\arg(B(x) v,v)|. $$
Then the operator $D_B = d + B^{-1} d^{*}_g B$ has a bounded
$S_\mu ^o$ holomorphic functional calculus in $\mH$. The constant
in this bound  depends only on $M$, $\|B\|$ and
$\kappa$.
\end{thm}

We begin the proof of Theorem \ref{thm:man1} with a localization
lemma. Let $\rho:U \longrightarrow B(0,4\delta)$ be a
diffeomorphism (or bi-Lipschitz mapping) for some open $U \subset
M$, $\delta>0$. Here we let $B(x,r)$ denote the ball in $\R^n$
with centre $x \in \R^n$ and radius $r>0$, where $n$ is the
dimension of $M$. Let $\rho^*$ denote the pullback by a function
$\rho$. Let $\Theta^B _t$ be as given in Definition \ref{defnope}
with $\Gamma:=d$ and $\Pi_B:=D_B$.

\begin{lem} \label{cov2}
We have $$ \int_0 ^1  \| \Theta^B _t u \|^2 \, \frac{dt}{t}
\lesssim \| u \|^2 $$ for every $u \in \mH$ with $\supp u \subset
\rho^{-1}( B(0,\delta))$. The bound here depends on $\delta$, the
hypothesis of Theorem \ref{thm:man1}, and the gradient bounds of
$\rho$ and $\rho^{-1}$.
\end{lem}

\begin{proof}
By Proposition \ref{pseudoloc} (adapted to the setting of a
compact Riemannian manifold) we have that
$$ \int_{M \setminus \rho^{-1}(B(0,2\delta)) } | \Theta^B _t u |^2 \, dx
\lesssim  t^2 \|u\|^2$$ and therefore that
$$
\int_0 ^1 \int_{M \setminus \rho^{-1}(B(0,2\delta)) } | \Theta^B
_t u |^2 \, dx \frac{dt}{t} \lesssim  \int_0 ^1  t^2 \| u \|^2 \,
\frac{dt}{t} \lesssim  \| u \|^2.$$ It remains to show that
\begin{equation} \label{eq:dtheta2} \int_0 ^1
\int_{\rho^{-1}(B(0,2\delta)) } | \Theta^B _t u |^2 \, dx
\frac{dt}{t} \lesssim  \| u \|^2 .\end{equation} We do this by
pushing the problem onto $\rnum^n$.

Let $\hat B$ be the multiplication operator on $L_2(\rnum^n ;
\wedge_\cnum \rnum^n)$ that coincides with the identity on
$\rnum^n \setminus B(0,4 \delta)$, and is otherwise fixed by the
condition that $(\rho^{-1})^* D_B \rho^* = D_{\hat B}$, where we
write $D_{\hat B}:= d+ (\hat B)^{-1} d^* \hat B$, and where $d^*$
denotes the adjoint of $d$ under the standard Euclidean metric.
Here $\hat B= (\rho_*/J_\rho)B\rho^*$, where
$\rho_*/J_\rho:L_2(\wedge T^*U)\rightarrow
L_2(B(0,4\delta);\wedge_\cnum\rnum^n)$ is the adjoint of
$\rho^*:L_2(B(0,4\delta);\wedge_\cnum\rnum^n)\rightarrow L_2(\wedge
T^*U)$ and $\rho_*$ denotes the pushforward and $J_\rho$ the
Jacobian determinant of $\rho$.
 By our hypotheses on $B$ we then have
$B_1 =(\hat B)^{-1}$ and $B_2=\hat B$ satisfy (H2,3,5) with bounds
that depend only on the hypotheses and the gradient bounds on
$\rho$ and $\rho^{-1}$. By Theorem \ref{mainthm} with $\{
\Gamma=d, (\hat B)^{-1}, \hat B\}$ we then have \begin{equation}
\label{missi} \int_0 ^1 \| t D_{\hat B} (\I + t^2 {D_{\hat
B}}^2)^{-1} v \|^2 \, \frac{dt}{t} \lesssim \| u \|^2
\end{equation} where, here and after we fix $v= (\rho^{-1})^* u$.

To complete the proof it suffices to show that
\begin{equation} \label{eq:diff1} \|(\rho^{-1} )^* it D_{B} ( \I +
t^2 {D_B}^2)^{-1} \rho^* v - it D_{\hat B} ( \I +t^2 {D_{\hat
B}}^2)^{-1} v \|_{L_2(B(0, 2\delta))} \lesssim  t  \| v
\|\end{equation} for every $0 < t \le 1$, and that these bounds
depend on the hypotheses and the gradient bounds on $\rho$ and
$\rho^{-1}$. (Indeed, we can then apply the triangle inequality
with (\ref{missi}) to the bound the left-hand side of
(\ref{eq:dtheta2}) by a controlled constant times $\|u\|^2 +
\int_0 ^1 t^2 \| u \|^2 \, \frac{dt}{t} \lesssim \| u \|^2.$) To
see (\ref{eq:diff1}) holds, let $\eta_1,\eta_2 : \rnum^n
\longrightarrow \rnum$ be smooth cut-off functions with
\begin{equation*}
\eta_i(x) = \left\{ \begin{array}{ll}
1    & \textrm{if $x \in B(0,(i+1) \delta)$}\\
0    & \textrm{if $x \in \rnum ^n \setminus B(0, (i+2) \delta)$} \\
\end{array} \right.
\end{equation*}
and $ | \nabla \eta_i | \le 2\delta^{-1}$ for $i=1,2$. Observe
that
\begin{equation*}
\begin{split} &(\I + it D_{\hat
B})^{-1} v(x) - (\rho^{-1})^*(\I+ it D_B )^{-1} \rho ^* v (x) \\
& \quad= (\rho^{-1})^*(\I + it D_B)^{-1} \rho^* \eta_2 \left(
(\rho^{-1})^* (\I + it D_B) \rho^* \eta_1 - (\I + it D_{\hat B} )
\right) ( \I + it D_{\hat B} )^{-1} v (x) \\
& \quad= (\rho^{-1})^*(\I + it D_B)^{-1} \rho^* \eta_2 (\I + it
D_{\hat B} )(\eta_1 -1 ) (\I+i t D_{\hat B} )^{-1} v (x) \\& \quad
= (\rho^{-1})^*(\I + it D_B)^{-1} \rho^* \eta_2 [itD_{\hat
B},\eta_1]  (\I+i t D_{\hat B} )^{-1} v (x)
\end{split}
\end{equation*}
for almost every   $x \in B(0,2\delta)$, and by Proposition
\ref{pseudoloc} has $L_2 (B(0,2\delta);\wedge_\C\rnum^n)$ norm
bounded by a constant multiple of $  t \| \nabla \eta_1  \| \|v \|
\lesssim  t \|v \|,$ where the constant depends only on the
assumed constants of the hypotheses. Estimate (\ref{eq:diff1}) now
follows by writing $Q^B _t= \tfrac{1}{2i} (-R^B _t + R^B _{-t} )$.
This completes the proof. \qedend \end{proof}

\begin{proof}[Proof of Theorem \ref{thm:man1}]
By Proposition \ref{whatisneeded}, we need to establish
(\ref{thesqest}) for every $u \in \ran(\Gamma)$, for each case
where $\{\Gamma,B^{-1},B\}$ is given by $\{d, B^{-1},B\}$,
$\{d^*_g, B,B^{-1}\}$, $\{d^*_g,B^*, (B^{-1})^*\}$ and $\{d,
(B^{-1})^{ *},B^*\}$. 
Let $H$ be the Hodge-star operator on $M$ and let $N$ be the operator
that changes sign of forms of odd degree. Then we have the unitary equivalence
$$H^*(d_g^* + B d B^{-1})H = Nd + \tilde B^{-1}(Nd)^* \tilde B$$
where $\tilde B= H^*B^{-1}H$ satisfies the same hypothesis as $B$.
Consequently, all four cases are essentially of the form
 $\{d, B^{-1},B\}$\ which we now consider.

Since $M$ is compact we can use Lemma \ref{cov2} with a
standard local chart/partition of unity argument to deduce that
\begin{equation*} \label{eq:dtheta1}
    \int_0 ^1 \|
\Theta^B _t u \|^2 \, \frac{dt}{t} \lesssim
\|u\|^2.\end{equation*} Again because $M$ is compact, and also
because $u \in \ran (d)$ and thus $P_t u \in \ran (D)$, we can
apply the Gaffney-G\aa rding inequality (see \cite[Theorem
7.3.2]{morrey:mult}) to deduce that $ \| P_t u \| \lesssim \| D
P_t u \|$, and therefore conclude that $$\int_1 ^\infty \|
\Theta^B _t P_t u \|^2 \, \frac{dt}{t} \lesssim  \int_1 ^\infty \|
      P_t   u \|^2 \, \frac{dt}{t} \lesssim \int_1 ^\infty \| t D P_t
u \|^2 \, \frac{dt}{t^3} \lesssim \int_1 ^\infty \| u \|^2 \,
\frac{dt}{t^3} \lesssim \|u\|^2.$$ This with Lemma \ref{cov2}
proves (\ref{thesqest}) and so completes the proof of Theorem
\ref{thm:man1}. \qedend \end{proof}

We now state an application of the above theorem. Given  a smooth
perturbation $g+h$ of $g$ we let
$$|h_x|=
\sup\{|h_x(v,v)| :  v\in \wedge T_x M \, ,\, g_x(v,v)=1\}$$ for
all $x \in M$, and define $\|h\|_\infty:=\sup_{x\in M}|h_x|$.
(This norm is equivalent to the one given in the Introduction, but
more useful for our purposes.)

\begin{thm} \label{thm:man2}
Let $M$ be a compact Riemannian manifold with metric $g$, let
$g+h$ be a measurable perturbation of $g$ with $\| h \|_\infty<
1/4$, and let $0 < \mu < \frac\pi2$ be given by $ \mu = \cos^{-1}
(1/4)$. Then we have
\begin{equation*}\|f (D_{g+h})-f(D_g)\| \lesssim \|f \|_\infty
\|h\|_\infty
\end{equation*}
for every bounded $f : S_\mu ^o \cup \{0 \} \longrightarrow \cnum$
holomorphic on $S_\mu ^o$. The constant in the above bound depends
only on $M$.
\end{thm}


\begin{rem} \label{deeper}
Lipschitz estimates like those in Theorem \ref{thm:man2} also hold
in terms of the quadratic estimates appearing in the second part
of Theorem \ref{prop:lips}. These  results are a consequence of
the deeper fact that the mapping given by $z \mapsto f (D_{g+zh})$
 depends holomorphically on $z \in \C$
when $|z|<\|h\|_\infty^{-1}$. These same results hold for any
manifold bi-Lipschitz equivalent  to Euclidean space, and follow
by arguments similar to those used in this section. We leave the
details to the reader.
\end{rem}

\begin{proof}[Proof of Theorem \ref{thm:man2}]
We can implicitly define   $A\in L_\infty(\mL_M)$ by the formula
$$ ((I+A(x))  u(x),  v(x) )_{g} = (u(x),v(x))_{g+h}$$
for every $u,v \in L_2(\wedge_\cnum T^* M)$. Here we let $(\cdot,
\cdot)_{g+h}$ and $(\cdot, \cdot)_{g}$ denote the metrics on $M$
corresponding to $g+h$ and $g$, respectively. Our hypothesis on
$g+h$ implies that $A\in L_\infty(\mL_M )$ with $ \| A \|_\infty =
\| h \|_\infty \le 1/4$ and therefore also
$$ \| I - (I+A)^{-1} \|_\infty \le \frac{\|h\|_\infty }{1 - \| h
\|_\infty} \le 1/3.$$
     Moreover, we have
$$ ((I+A)   d^*_{g+h} u,   v )_{g} = (d^*_{g+h} u , v)_{g+h} =
(u,dv)_{g+h}  = ((I+A)   u, d   v)_{g} = ( d^*_{g} (I+A)   u,
v)_{g} $$ for every $u,v \in L_2(\wedge_\cnum T^* M)$ with $u \in
\dom (d^*_{g+h})$ and $v \in \dom(d)$, and therefore
$$  D_{g+h} =   d  +
d^*_{g+h} =  d +  (I+A)^{-1} d^*_{g} (I+A) .$$ The  desired
result now follows from an application of Theorem \ref{thm:man1}
and results analogous to Theorem \ref{prop:lips} with $A_2 = A$,
$A_1 = (I +A)^{-1}-I$, $\eta_i=1/2$, $\kappa_i=1$, and
$B_i=I$ for $i=1,2$. \qedend \end{proof}

\begin{appendix}

\section{Further properties of the Hodge decomposition} \label{sect:proj}
\label{append}

In this appendix we verify the claim of Remark \ref{remobs}. As in
Section \ref{opthsec}, we assume that the triple of operators
$\{\Gamma, B_1,B_2 \}$ in a Hilbert space $\mH$ satisfies
properties (H1--3). We begin with a lemma.

\begin{lem}\label{lem:multi}
The Hodge projections can be represented as limits of resolvents
in the following ways:
$$\PP^0 u  = \lim_{n \to \infty} (\I + i n \Pi_B)^{-1} u = \lim_{n
\to \infty} (\I - i n \Pi_B)^{-1}u\ \quad \text{for all} \quad
u\in\mH\ ;$$ $$ \PP^1 u  = \lim_{n \to \infty} in\Gamma^* _B(\I +
in\Pi_B)^{-1}u = \lim_{n \to \infty} i n \Gamma_B ^* (-\I + i n
\Pi_B)^{-1}u\ \quad  \text{for all} \quad u\in\mH \ ;$$ $$ \PP^1 u
= \lim_{n \to \infty} (\I + i n \Pi_B)^{-1} in \Gamma u= \lim_{n
\to \infty} (-\I+ i n \Pi_B)^{-1} i n \Gamma u\ \quad \text{for
all} \quad u\in\dom(\Gamma)\ ;$$ $$ \PP^2u  = \lim_{n \to \infty}
in\Gamma (\I + in\Pi_B)^{-1} u= \lim_{n \to \infty} i n \Gamma
(-\I + i n \Pi_B)^{-1}u\ \quad \text{for all} \quad u\in\mH \ ;$$
$$ \PP^2 u = \lim_{n \to \infty} (\I + i n \Pi_B)^{-1} in
\Gamma_B^*u = \lim_{n \to \infty} (-\I+ i n \Pi_B)^{-1} i n
\Gamma_B^*u\ \quad \text{for all} \quad u\in\dom(\Gamma^*_B)\ .$$
\end{lem}

\begin{rem}

If $\nul(\Pi_B)=\{0\}$, then
$\PP^1=\Gamma^*_B\Pi_B^{-1}=\Pi_B^{-1}\Gamma$ and
$\PP^2=\Gamma\Pi_B^{-1}=\Pi_B^{-1}\Gamma^*_B$ on the appropriate
domains, in which case the proofs would be somewhat more direct.
\end{rem}

Note that, by Proposition  \ref{typeomega} and Lemma
\ref{lem:gam1}, each of the operator sequences $(\I + i n
\Pi_B)^{-1}$, $in\Gamma^* _B(\I + in\Pi_B)^{-1}$, etc, is
uniformly bounded in $n$. It is not a-priori clear that their
strong limits exist. This will be shown in the course of the
proof.

\begin{proof}
We begin by showing that
\begin{equation} \label{eq:strong1} Q_n^B u = n \Pi_B (\I
       +n^2 {\Pi_B}^2)^{-1} u \to 0\end{equation} as $ n \to
       \infty$ for every $u \in \mH$. The expression on the left vanishes
       if $u \in  \nul(\Pi_B)$, so it suffices to consider the case
when $u \in \overline{\ran(\Pi_B)}$. If $u=\Pi_B v \in
\ran(\Pi_B)$, then $$ \|  Q_n^B u\| =\|Q_n^B\Pi_B v\| =
\tfrac1n\|v-P_n^Bv\| \lesssim \tfrac1n\|v\|\to 0$$ as $n \to
\infty$. Since, by Proposition \ref{typeomega}, the sequence
        $\| Q_n^B \|$
is uniformly bounded, we conclude by a standard continuity
argument that (\ref{eq:strong1}) holds for every $u \in
\overline{\ran(\Pi_B)}$.

Define operators $T_0, T_1$ and $T_2$ on $\mH$ by
$$T_0 u= \lim_{n \to \infty} (\I+in \Pi_B)^{-1}u ,$$
$$ T_1 u = \lim_{n \to \infty} in \Gamma^*_B (\I + in \Pi_B)^{-1} u
\quad \text{and} \quad T_2 u= \lim_{n \to \infty} in \Gamma (\I +
in \Pi_B)^{-1}u$$ whenever $u \in \mH$ and the corresponding limit
exists. We next  show that
\begin{equation} \label{eq:redu11} T_0u = \lim_{n \to \infty} (\I -
i n \Pi_B)^{-1}u,\end{equation}
\begin{equation} \label{eq:redu13}  \begin{split}
\PPT_1 u = \lim_{n \to \infty} i n \Gamma_B ^* (-\I + i n
\Pi_B)^{-1} u &= \lim_{n \to \infty} (\I + i n \Pi_B)^{-1} in
\Gamma u
       \\ &= \lim_{n \to \infty} (-\I+ i n \Pi_B)^{-1} i n \Gamma u
\end{split}
\end{equation}
       and
\begin{equation} \label{eq:redu12}  \begin{split}
\PPT_2 u = \lim_{n \to \infty} i n \Gamma (-\I + i n \Pi_B)^{-1} u
&= \lim_{n \to \infty} (\I + i n \Pi_B)^{-1} in \Gamma_B^* u\\ &=
\lim_{n \to \infty} (-\I+ i n \Pi_B)^{-1} i n \Gamma_B^* u
\end{split}
\end{equation}
whenever $u \in \mH$ (and when required, $u\in \dom(\Gamma)$ or
$u\in\dom(\Gamma^*_B)$)
       and the corresponding limit exists. Here we
interpret the above as saying that if one such limit exists, then
the limits that are indicated to be equal,  also exist.

Equation (\ref{eq:redu11}) follows by (\ref{eq:strong1}) and the
fact that
\begin{equation*} \label{eq:vani1} (\I+ i n \Pi_B)^{-1} -
(\I - i n \Pi _B)^{-1}  = -2in \Pi_B ( \I + n^2 \Pi _B ^2 )^{-1}\
.\end{equation*}

       To see the first
equality in (\ref{eq:redu12}), observe that by (\ref{eq:strong1})
and Lemma \ref{lem:gam1} we have
\begin{equation*}
\begin{split} \| in\Gamma (\I + i n \Pi_B)^{-1}u - in \Gamma(-\I + in
\Pi_B)^{-1} u\| &= \| 2 in \Gamma (\I+ n^2 {\Pi_B}^2)^{-1}u \| \\&
\lesssim \, \| n \Pi_B (\I+ n^2 {\Pi_B}^2)^{-1}  u \| \to 0
\end{split}
\end{equation*} as $n \to \infty$.  The second equality in (\ref{eq:redu12})
follows from (\ref{eq:strong1}) and the identity
\begin{equation} \label{eq:gamcom}
\begin{split}  &i n \Gamma(-\I+ in \Pi_B)^{-1}u - (\I+in\Pi_B)^{-1} in
\Gamma_B ^* u \\ &= ( \I+ in \Pi_B)^{-1} \left( (\I+ in \Pi_B) in
\Gamma - in \Gamma_B ^*  (-\I + in \Pi_B) \right) (-\I + in
\Pi_B)^{-1} u \\ &= ( \I+ in \Pi_B)^{-1}
       (in \Gamma + in \Gamma_B ^*) (-\I + in \Pi_B)^{-1} u \\&= -in \Pi_B
       (\I
       +n^2 {\Pi_B}^2)^{-1} u =  -i Q_n^B u
\end{split}
\end{equation}
for all $ u\in \dom(\Gamma^*_B) $.

The remaining equality in (\ref{eq:redu12}) as well as Equation
(\ref{eq:redu13}) can be proved by similar arguments.

We note that $T_0 u = u $ when $u\in\nul(\Pi_B)$, and, by adapting
the proof of (\ref{eq:strong1}), that $T_0 u = 0 $ when
$u\in\ran(\Pi_B)$ and hence when $u\in\overline{\ran(\Pi_B)}$.
Therefore $T_0=\PP^0$.

Now investigate $T_1$. By (\ref{eq:redu13}), $T_1u=0$ when
$u\in\nul(\Gamma)$.  If $u\in\ran(\Gamma^*_B)$, let $u=\Gamma^*_B
v$, where, by Proposition~\ref{hodgedec}, we may assume that
$v\in\clos{\ran(\Gamma)}$. Using the facts that $T_0 v=0$ and that
$\Gamma^*_B$ is closed, we obtain
\begin{equation}  T_1 u= \lim_{n\rightarrow\infty}
in\Gamma^*_B(\I+in\Pi_B)^{-1}\Pi_B v
        = \lim_{n\rightarrow\infty} \Gamma^*_B(\I-(\I+in\Pi_B)^{-1}) v
        = \Gamma^*_B v=u\,.
\end{equation}
By a standard argument, we find that $T_1 u = u$ when
$u\in\overline{\ran(\Gamma^*_B)}$. Therefore $T_1=\PP^1$.

Similarly, $T_2 u = 0 $ when $u\in\nul(\Gamma^*_B)$, and $T_2 u =
u $ when $u\in\overline{\ran(\Gamma)}$, so that $T_2=\PP^2$.
\qedend
\end{proof}

Define operators $\Pa$ and $\Pb$ on $\mH$ by
$$\Pa u  = \lim_{n \to \infty} in\Gamma^*B_2(\I + in\Pi_B)^{-1}u
\quad  \text{for all} \quad u \in \mH\,,$$
$$\Pb   u= \lim_{n \to
\infty} (\I + i n \Pi_B)^{-1} in B_1\Gamma^* u \quad \text{for
all}  \quad u \in \dom(\Gamma^*)\,.$$ The fact that the limits
defining $\Pa$ and $\Pb$ exist and define
   bounded operators, as well as the fact that
$\PP^1=B_1\Pa$ and $\PP^2=\Pb B_2$, now follow
  from (\ref{equiv1}),
(\ref{equiv2}). We remark that for (\ref{eq:deri1}),
(\ref{eq:deri2}) and (\ref{eq:deri3}) to be true, the sum of the
right-hand sides must equal zero, which requires
$$\PP^2A_1\Pa+\Pb A_2\PP^1=0.$$ Indeed, this is a consequence of
the assumption $\Gamma^*B_2B_1\Gamma^*=0$.

\begin{proof}[Proof of Remark \ref{remobs}]

Let $T_0^n=(\I+in\Pi_B)^{-1}$, $T_1^n= in \Gamma^*_B(\I+ in
\Pi_B)^{-1} $ and $T_2^n= in\Gamma(\I+ in \Pi_B)^{-1}$. By
Proposition~\ref{hodgedec} and Lemma~\ref{lem:gam1} we have that
the mappings $z \mapsto T_0^n$, $z \mapsto T_1^n$ and $z \mapsto
T_2^n$ are uniformly bounded. Furthermore Lemma \ref{lem:multi}
shows that $T_0^n\rightarrow \PP^0$, $T_1^n\rightarrow \PP^1$ and
$T_2^n\rightarrow \PP^2$ strongly. Thus it will follow from Lemma
\ref{lem:oper} that $z \mapsto \mathbf{P}^0 _{B(z)}$, $z \mapsto
\mathbf{P}^1 _{B(z)}$ and $z \mapsto \mathbf{P}^2 _{B(z)}$ are
holomorphic with derivatives as stated in (\ref{eq:deri1}),
(\ref{eq:deri2}) and (\ref{eq:deri3}) once we prove that $T_i^n$,
$i=1,2,3$ are holomorphic functions and that $\frac{d}{dz} T_i^n$
have as strong limits the right hand sides in (\ref{eq:deri1}),
(\ref{eq:deri2}) and (\ref{eq:deri3}) respectively.

For $T_0^n$, we see that \begin{equation} \label{eq:deri8}
\begin{split} \frac{d}{dz} T^n_0u=\frac{d}{dz} (\I + in \Pi_B)^{-1}u
&= -(\I + in \Pi_B)^{-1} A_1 in \Gamma^* B_2 (\I + in \Pi_B)^{-1}u
\\ &\quad- (\I + in \Pi_B)^{-1}  B_1in \Gamma^* A_2 (\I + in
\Pi_B)^{-1}u \\
&\to (-\PP ^0 A_1 \Pa - \Pb  A_2 \PP^0)u\,.
\end{split}
\end{equation}

For $T_1^n$, we see that when $u\in\dom(\Gamma)$,
\begin{equation} \label{eq:deri4} \begin{split}
\frac{d}{dz}T^n_1u=\frac{d}{dz} (\I+ in \Pi_B)^{-1} in \Gamma u &=
-(\I + in \Pi_B)^{-1} A_1 in \Gamma^* B_2 (\I + in \Pi_B)^{-1} in
\Gamma u\\ &\quad- (\I + in \Pi_B)^{-1} B_1in \Gamma^* A_2 (\I +
in \Pi_B)^{-1} in \Gamma u\,.
\end{split}
\end{equation}
The second term on the right-hand side converges to $-\Pb A_2
\PP^1 u$. In order to calculate the first term on the right-hand
side, we note by an argument similar to (\ref{eq:gamcom}) that
\begin{equation*}
\begin{split}  &i n \Gamma^*_B(\I+ in \Pi_B)^{-1} - (-\I+in\Pi_B)^{-1} in
\Gamma = -in \Pi_B
       ( \I+ in \Pi_B)^{-1} (-\I + in \Pi_B)^{-1} . \end{split} \end{equation*}
This with the fact that $\Gamma^2=0$ implies that
$$in \Gamma^* _B (\I + in \Pi_B)^{-1} in \Gamma =-in \Pi_B
       ( \I+ in \Pi_B)^{-1} (-\I + in \Pi_B)^{-1}  in \Gamma
        \to ( \PP^0 -I) \PP^1 = -\PP^1
$$
as $n \to \infty$. Therefore the  first term on the right-hand
side of (\ref{eq:deri4}) converges to $\PP^0 A_1 \Pa u$ as $n \to
\infty$.

A similar argument shows that $
   \frac{d}{dz}T_2^n u\to (-\PP ^2 A_1 \Pa + \Pb A_2 \PP^0)u.
$ This completes the proof. \qedend
\end{proof}

\end{appendix}

\bibliographystyle{acm}

\end{document}